\documentclass[letterpaper,11pt]{article}

\usepackage{amsfonts}
\usepackage{amsthm}
\usepackage{amssymb}
\usepackage{amsmath}
\usepackage{mathrsfs}
\usepackage{verbatim}
\usepackage{bm}
\usepackage{graphicx}
\usepackage{color}
\usepackage[left=1.0in,top=1.0in,bottom=1.0in,right=1.0in]{geometry}

\newtheorem{theorem}{Theorem}[section]
\newtheorem{assumption}[theorem]{Assumption}

\newtheorem{definition}[theorem]{Definition}
\newtheorem{lemma}[theorem]{Lemma}
\newtheorem{proposition}[theorem]{Proposition}
\newtheorem{remark}[theorem]{Remark}

\numberwithin{equation}{section}

\definecolor{Red}{rgb}{1.00, 0.00, 0.00}

\definecolor{Blue}{rgb}{0,0,1}

\definecolor{Green}{rgb}{0.2, 0.5, 0.2}

\def\Swiech{{\accent"13S}wie{\hbox{\kern -0.21em\lower 0.79ex\hbox{$\textfont1=\scriptfont1\lhook$}}}ch}

\title{A Viscosity Approach to a Stochastic Control Problem on a Bounded Domain}

\author{Ruoting Gong\thanks{Department of Applied Mathematics, Illinois Institute of Technology, Chicago, IL 60616, U.S.A. {\tt rgong2@iit.edu}} \and Christian Houdr\'{e}\thanks{School of Mathematics, Georgia Institute of Technology, Atlanta, GA 30332, U.S.A. {\tt houdre@math.gatech.edu}. Research supported in part by the grants \# 246283 and \# 524678 from the Simons Foundation.}}

\date{October 21, 2017}

\begin{document}

\maketitle

\begin{abstract}

We study a stochastic control problem on a bounded domain, which arises from a continuous-time optimal management model. Via the corresponding Hamilton-Jacobi-Bellman equation the value function is shown to be jointly continuous and to satisfy the Dynamic Programming Principle. These properties directly lead to the conclusion that the value function is a viscosity solution to the Hamilton-Jacobi-Bellman equation. Uniqueness of the solution is then also established.

\end{abstract}

\vspace{0.3cm}


\vspace{0.3cm}


\section{Introduction}

Stochastic representation formulas establish natural connections between the study of stochastic processes, and partial differential equations (PDEs). Most notably, the dynamic programming principle and the theory of regular and viscosity solutions establish rigorous connection between stochastic optimal control problems and fully nonlinear Hamilton-Jacobi-Bellman (HJB) equations. Thus, value functions of optimal control problems are identified as regular or viscosity solutions to the associated HJB equations and, in particular, provide stochastic representations to those solutions. Such techniques found applications in many areas, such as finance, economics, physics, biology, and engineering. Various results on regular and viscosity solutions to HJB and Isaacs PDEs in bounded and unbounded domains and their associated stochastic control problems can be found, for instance, in~\cite{BuckdahnLi:2008}, \cite{FlemingSoner:2006}, \cite{FlemingSouganidis:1989}, \cite{Katsoulakis:1995}, \cite{Kovats:2009}, \cite{Krylov:1980}, \cite{Lions:1983(1)}, \cite{Lions:1983(3)}, \cite{MaZhang:2002}, \cite{Nisio:2015}, \cite{PardouxPeng:1992}, \cite{Swiech:1996}, \cite{Touzi:2002}, \cite{Touzi:2013}, \cite{YongZhou:1999}, \cite{Zhang:2005}, and the references therein. The literature on the subject is huge.

In this paper, we consider a stochastic control problem on a bounded domain, arising from an optimal management model. We would like to emphasize that, although some components of our model (e.g., the value process of the underlying project) are similar to those of~\cite{GiatSubramanian:2013} in which a dynamic principal-agent model is introduced and investigated, our optimal control problem is different from the classical principal-agent problem as studied, for example, in~\cite{CvitanicPossamaiTouzi:2017}, \cite{CvitanicZhang:2013} and~\cite{LaffontMartimort:2002}. In particular, there is no dynamic contracting structure in our model, and the manager takes charge of all control variables. Our main focus is to identify the value function of our stochastic control problem as the unique viscosity solution to the Dirichlet terminal-boundary value problem for the associated degenerate HJB equation in a bounded domain. This is a classical problem which is very technical and whose full details are often omitted or overlooked, especially for problems in bounded domains. Our method is similar to that of~\cite{FlemingSoner:2006}. We establish the joint continuity of the value function as well as the dynamic programming principle, from which the value function can be directly verified to be a viscosity solution to the associated HJB equation. We also establish the comparison principle for the HJB equation using the well-known Ishii's lemma. The main difficulties come from the fact that we are dealing with a degenerate HJB equation on a bounded region. The degeneracy of the HJB equation is tackled by approximating our PDE by equations which are non-degenerate, have more regular coefficients, and are considered on slightly smaller domains with smooth boundaries. Such equations have classical solutions which can be identified as value functions of associated optimal control problems. We then pass to the limits with various approximations. The bounded region forces a lot of technical estimates involving the analysis of the behavior of stochastic processes and their exit times. We remark that making the HJB equation non-degenerate by adding a small Laplacian term to the equation corresponds to the introduction of another independent Wiener process on the level of the stochastic control problem, and hence to possible enlargement of the reference probability space.

The present paper is organized as follows. Section \ref{sec:Survey} provides a brief review of the literature on stochastic control and viscosity solutions. In Section \ref{sec:Model}, we state the optimal management model, and formulate the stochastic control problem and the corresponding HJB equation. Section \ref{sec:Existence} and Section \ref{sec:Uniqueness} contain the main results of the manuscript. In Section \ref{sec:Existence}, we first prove the joint continuity of the value function and establish the Dynamic Programming Principle. We then verify that the value function is indeed a viscosity solution of the HJB equation with terminal/boundary condition. Finally, in Section \ref{sec:Uniqueness}, we establish the uniqueness of the viscosity solution of the HJB equation under polynomial growth.

\medskip
\noindent
\textbf{Acknowledgement:} It is a pleasure to thank A. Subramanian who introduced us to~\cite{GiatSubramanian:2013} and A. Swiech for discussions, bibliographical help, and setting us straight.

\section{A Brief Literature Survey}\label{sec:Survey}

Our study is mainly concerned with stochastic control and viscosity solutions. Many authors have introduced different notions of generalized solutions in order to prove that the value function is a solution of the corresponding HJB equation. Kr\v{u}zkov~\cite{Kruzkov:1960}$-$\cite{Kruzkov:1970} built a systematic theory for first-order Hamilton-Jacobi (HJ) equations with smooth and convex Hamiltonians. Fleming~\cite{Fleming:1964}-\cite{Fleming:1969} independently introduced the vanishing of viscosity, combined with the differential games techniques, to study the HJ equation. Clarke and Vinter~\cite{ClarkVinter:1983} used Clarke's notion of generalized gradients to introduce generalized solutions of the HJB equations. In that framework, the HJB equation can have more than one solution, and the value function is one of them. However, generalized gradients may not be readily used to solve second-order HJB equations corresponding to stochastic problems. A survey of HJB equations is presented in~\cite{BardiCapuzzoDolcetta:1997}.

In the early 1980s, Crandall and Lions~\cite{CrandallLions:1983} introduced the notion of viscosity solution for first-order HJB equations. The first treatment of viscosity solutions of second-order dynamic programming equations was then given by Lions~\cite{Lions:1983(1)}-\cite{Lions:1983(3)}, who investigated the degenerate second-order HJB equation using a Feynman-Kac-type technique, and represented solutions as value functions of some stochastic control problems. For general second-order equations which are not necessarily dynamic programming equations, this technique is no longer appropriate. Jensen~\cite{Jensen:1988} first proved a uniqueness result for a general second-order equation, in which semiconvex and concave approximations of a function are given by using the distance to the graph of this function. Another important step in the development of the second-order problems is Ishii's Lemma (cf.~\cite{Ishii:1989}). Since then the theory of second-order equations has seen great progress. In particular, the analytical results of Crandall and Ishii~\cite{CrandallIshii:1990} have been used in almost all comparison results. We refer to the survey article~\cite{CrandallIshiiLions:1992} for more detailed information. Fleming and Soner~\cite{FlemingSoner:2006} provided a rigorous approach to the control theory of Markov diffusion processes. Specifically, when uniform ellipticity is satisfied, the value function is shown to be a classical solution of the corresponding second-order HJB equation. When uniform ellipticity is abandoned, a systematic analysis of the value function is provided and a strong version of the dynamic programming principle is established, reducing the problem to the uniformly elliptic case. Similar results are summarized and developed in~\cite{YongZhou:1999} via independent approaches. Viscosity solutions to HJB integro-PDEs and their stochastic representation formulas as value functions of the associated stochastic optimal control problems were originally investigated by Soner~\cite{Soner:1986}, \cite{Soner:1988}. Stochastic optimal control of jump-diffusion processes and various results on the associated HJB equations are discussed in~\cite{OksendalSulem:2007}.

\section{Basic Settings and Preliminary Results}\label{sec:Model}

\subsection{The Optimal Management Model}\label{subsec:PrinAgModel}

Throughout this paper, let $T>0$ be a fixed terminal time. In our model, a risk-neutral manager with capital carries out a project with a group of collaborators on the finite horizon $[0,T]$; the manager and his group being referred to as ``the team". The project can potentially generate value through capital investments from the manager and human effort investments from the team. The key variable in the model is the project's terminal value process $V:=(V_{t})_{t\in[0,T]}$. The incremental termination value $dV_{t}$, i.e., the change in the termination value over the infinitesimal period $[t,t+dt]$, is the sum of a base output which is unaffected by the actions of the team, and a discretionary output which depends on the manager's capital investments and the team's efforts.

More precisely, let $(\Omega,\mathscr{F},\mathbb{P})$ be a complete probability space, where $\mathbb{P}$ is the market probability measure. Let $W:=(W_{t})_{t\in[0,T]}$ be a standard Brownian motion on $(\Omega,\mathscr{F},\mathbb{P})$, and let $\mathbb{F}_{T}:=(\mathscr{F}_{t})_{t\in[0,T]}$ be the augmented filtration generated by $W$. Let $\Xi:=(\Xi_{t})_{t\in[0,T]}$, $L:=(L_{t})_{t\in[0,T]}$, and $C:=(C_{t})_{t\in[0,T]}$ be $\mathbb{F}_{T}$-progressively measurable processes, describing respectively the manager's investments, the team's efforts, and the project's consumptions over time. Throughout this article, we assume that $(\Xi,L,C)$ takes value in a compact set $\mathscr{U}:=[0,\mathscr{N}]\times[0,\mathscr{L}]\times[0,\mathscr{C}]$, where $\mathscr{N}$, $\mathscr{L}$, and $\mathscr{C}$ are positive constants. Intuitively, $\mathscr{N}$ represents the maximal capital investment the manager can afford, $\mathscr{L}$ represents the maximal effort the team can make, and $\mathscr{C}$ represents the maximal possible consumptions of the project.

The following construction of the value process $V$ is similar to~\cite{GiatSubramanian:2013}. Let $\Theta:=(\Theta_{t})_{t\in[0,T]}$ be a diffusion process representing the intrinsic quality of the underlying project. Under $\mathbb{P}$, $\Theta$ is assumed to evolve as:
\begin{align*}
d\,\Theta_{t}=\vartheta(t)\,dt+\sigma(t)\,dW_{t},
\end{align*}
where $\vartheta$ and $\sigma$ are respectively the deterministic drift and volatility function, with $\sigma(t)>0$ for all $t\in[0,T]$. The value process of the project, under $\mathbb{P}$, is then defined as
\begin{align*}
V_{t}=\varrho\,W(t),\quad t\in[0,T],
\end{align*}
where $\varrho^{2}$ is the intrinsic risk of the project with $\varrho>0$. Respectively define $Z:=(Z_{t})_{t\in[0,T]}$ and $B:=(B_{t})_{t\in[0,T]}$ via:
\begin{align}\label{eq:Girsanov}
Z_{t}&:=\exp\left(\int_{0}^{t}\!\left(\Theta_{s}\!+\!A\,\Xi_{s}^{\alpha}L_{s}^{\beta}\!-\!C_{s}\right)\!\varrho^{-1}dW_{s}-\frac{1}{2}\int_{0}^{t}\!\left(\Theta_{s}\!+\!A\,\Xi_{s}^{\alpha}L_{s}^{\beta}\!-\!C_{s}\right)^{2}\!\!\varrho^{-2}ds\right),\quad t\in[0,T],\\
B_{t}&:=W_{t}-\int_{0}^{t}\left(\Theta_{s}+A\,\Xi_{s}^{\alpha}L_{s}^{\beta}-C_{s}\right)\varrho^{-1}\,ds,\quad t\in[0,T],\nonumber
\end{align}
where $\Phi(\xi,\ell):=A\xi^{\alpha}\ell^{\beta}$ is the Cobb-Douglas production function (cf.~\cite{CobbDouglas:1928}) with $\alpha>0$ and $\beta>0$. By Girsanov's Theorem\footnote{We will add requirements on $\Theta$ so that the Novikov condition holds true, see Remark \ref{rem:Girsanov} below.}, under the new probability measure $\mathbb{Q}$, with
\begin{align}\label{eq:ProbMeasChange}
\frac{\left.d\mathbb{Q}\right|_{\mathscr{F}_{t}}}{\left.d\mathbb{P}\right|_{\mathscr{F}_{t}}}=Z_{t},\quad t\in[0,T],
\end{align}
$B$ is a standard Brownian motion. The new probability measure $\mathbb{Q}$ represents the manager's belief towards the market. Moreover, under $\mathbb{Q}$, the value process $V$ evolves as:
\begin{align}
dV_{t}&=\left(\Theta_{t}+A\,\Xi_{t}^{\alpha}L_{t}^{\beta}-C_{t}\right)dt+\varrho dB(t)\nonumber\\
\label{eq:ValProc} &=\underbrace{\Theta_{t}\,dt+\varrho\,dB_{t}}_{\text{base output}}+\underbrace{\left(\Phi(\Xi_{t},L_{t})-C_{t}\right)dt}_{\text{discretionary output}}.
\end{align}
As in~\cite{GiatSubramanian:2013}, the value process $V$ is only accessible to the manager under $\mathbb{Q}$ in the form of \eqref{eq:ValProc}.

Let $P:=(P_{t})_{t\in[0,T]}$ be the payoff process to the manager and the collaborators. Assume that the team has a minimal payoff tolerance $R>0$. The project is feasible at time $t\in[0,T]$ if the team is guaranteed the minimal payoff $R$ at time $t$, i.e., $P_{t}>R$. In the sequel, we will model $P$ as a diffusion process under $\mathbb{Q}$ (and thus under $\mathbb{P}$) whose drift and volatility depend on the manager's investment, the efforts of the team, and the consumptions of the project. Moreover, the project is assumed to incur a disutility of the team's effort. The rate of disutility from the team's effort in the period $[t,t+dt]$ is modeled as $\kappa\,L_{t}^{\gamma}\,dt$ with $\kappa>0$ and $\gamma>0$. Hence, if the project is terminated at time $\tau$, where $\tau$ is an $(\mathscr{F}_{t})_{t\in[0,T]}$-stopping time (see \eqref{eq:DefTau} for the exact definition of the random terminal time), then the manager's expected utility is given by
\begin{align}\label{eq:ExpUtl}
\mathbb{E}_{\mathbb{Q}}\left(\int_{0}^{\tau}\kappa L_{t}^{\gamma}\,dt-P_{\tau}\right)=\mathbb{E}_{\mathbb{P}}\left(\int_{0}^{\tau}\kappa L_{t}^{\gamma}Z_{t}\,dt-P_{\tau}Z_{\tau}\right).
\end{align}
Our goal is now to find an optimal triplet $(\Xi,L,C)$ which minimizes \eqref{eq:ExpUtl}.

\subsection{The Stochastic Control Problem}\label{subsec:StoContProb}

Let us first reiterate our setting. Let $(\Omega,\mathscr{F},(\mathscr{F}_{t})_{t\in[0,T]},\mathbb{P})$ be a complete filtered probability space, on which a standard Brownian motion $W=(W_{t})_{t\in[0,T]}$ is defined. Here $(\mathscr{F}_{t})_{t\in[0,T]}$ is the complete augmented filtration generated by $W$, and it thus satisfies the usual conditions. Let $\mathscr{U}=[0,\mathscr{N}]\times[0,\mathscr{L}]\times[0,\mathscr{C}]$ be the control domain with $\mathscr{N}>0$, $\mathscr{L}>0$ and $\mathscr{C}>0$, and let
\begin{align*}
\mathcal{U}[0,T]:=\left\{U=(\Xi,L,C):[0,T]\times\Omega\rightarrow\mathscr{U};\,\,(U_{t})_{t\in[0,T]}\,\,\text{is }(\mathscr{F}_{t})_{t\in[0,T]}-\text{progressively measurable}\right\}.
\end{align*}
Let $\mathscr{O}:=(R,\infty)\times(0,\infty)\times(-H,H)$ with $R>0$ and $H>0$. On $\mathscr{O}$, we consider the following SDE system for $X:=(X_{t})_{t\in[0,T]}$, where $X_{t}:=(P_{t},Z_{t},\Theta_{t})^{T}$, $t\in[0,T]$, under the control process $U=(\Xi,L,C)\in\mathcal{U}[0,T]$:
\begin{align}\label{eq:StateEqP}
dP_{t}&=b\left(t,P_{t},\Xi_{t},L_{t},C_{t}\right)dt+\Sigma\left(t,P_{t},\Xi_{t},L_{t},C_{t}\right)dW_{t}\\
\label{eq:StateEqZ} dZ_{t}&=-\varrho^{-1}Z_{t}\left(\Theta_{t}+A\,\Xi_{t}^{\alpha}L_{t}^{\beta}-C_{t}\right)dW_{t}\\
\label{eq:StateEqTheta} d\Theta_{t}&=\vartheta(t)\,\zeta_{H}(\Theta_{t})\,dt+\sigma(t)\,\zeta_{H}(\Theta_{t})\,dW_{t},
\end{align}
with initial condition
\begin{align}\label{eq:StateEqInit}
X_{0}=x:=(p,z,\theta)\in\overline{\mathscr{O}},
\end{align}
where $b,\Sigma:[0,T]\times\mathbb{R}\times\mathscr{U}\rightarrow\mathbb{R}$ and $\vartheta,\sigma:[0,T]\times\mathbb{R}$ are Borel-measurable functions, where $\varrho$, $A$, $\alpha$, and $\beta$ are positive constants, and where, as usual, $\overline{\mathscr{O}}$ is the closure of $\mathscr{O}$. Moreover, $\zeta_{H}:\mathbb{R}\rightarrow[0,1]$ is a deterministic $C^{2}$ function such that
\begin{align*}
\zeta_{H}(\theta)=0,\,\,\,\text{for }\,|\theta|\geq H+1;\quad\zeta_{H}(\theta)=1,\,\,\,\text{for }\,|\theta|\leq H.
\end{align*}
\begin{remark}\label{rem:Girsanov}
The choice of $\zeta$ above ensures the boundedness of the process $\Theta$, which, in turn, ensures the validity of the Novikov condition for $Z$, so that the Girsanov change of measure \eqref{eq:ProbMeasChange} is valid.
\end{remark}
In what follows, we will use $\mathbb{P}^{t,x}$ and $\mathbb{E}^{t,x}$ to respectively denote the probability and the expectation with respect to the initial data $x\in\overline{\mathscr{O}}$ starting at time $t\in[0,T]$. When $t=0$, we will omit the superscript $t$. For simplicity of notations, with $x=(p,z,\theta)\in\overline{\mathscr{O}}$ and $u=(\xi,\ell,c)\in\mathscr{U}$, let now
\begin{align*}
\vec{f}(s,x,u)&:=\left(b(s,p,u),\,0,\,\vartheta(s)\,\zeta_{H}(\theta)\right)^{T},\\
\vec{\sigma}(s,x,u)&:=\left(\Sigma(s,p,u),-\varrho^{-1}z\left(\theta+A\xi^{\alpha}\ell^{\beta}-c\right),\sigma(t)\zeta_{H}(\theta)\right)^{T},\\
\vec{a}(s,x,u)&:=\vec{\sigma}(s,x,u)\,\vec{\sigma}^{T}(s,x,u),\\
\mathcal{L}(x,u)&:=\kappa\ell^{\gamma}z.
\end{align*}
Using the above notations, the state equations \eqref{eq:StateEqP}-\eqref{eq:StateEqTheta} can be rewritten as
\begin{align}\label{eq:StateEq}
dX_{t}=\vec{f}(t,X_{t},U_{t})\,dt+\vec{\sigma}(t,X_{t},U_{t})\,dW_{t}.
\end{align}
For any $U\in\mathcal{U}[0,T]$ and any $x\in\overline{\mathscr{O}}$, define the cost function
\begin{align}\label{eq:CostFunt}
J(x;U):=\mathbb{E}^{x}\left(\int_{0}^{\tau}\mathcal{L}(X_{t},U_{t})\,dt-P_{\tau}Z_{\tau}\right),
\end{align}
where $\kappa>0$ and $\gamma>0$ are constants, and where
\begin{align}\label{eq:DefTau}
\tau=\tau(x;U):=\inf\left\{t\in[0,T]:\,X_{t}\not\in\mathscr{O},\,X_{0}=x\right\}\wedge T,
\end{align}
with the convention $\inf(\emptyset)=+\infty$.

\medskip
\noindent
\textbf{Problem (SC).}$\,\,\,$ Minimize \eqref{eq:CostFunt} over $\mathcal{U}[0,T]$. That is, for each fixed $x\in\overline{\mathscr{O}}$, find $U^{\ast}=U^{\ast}(x)\in\mathcal{U}[0,T]$, such that
\begin{align*}
J(x;U^{\ast})=\inf_{U\in\mathcal{U}[0,T]}J(x;U).
\end{align*}
\begin{remark}\label{rem:StoContWeakForm}
To search for the optimal control, we consider the weak formulation of the above stochastic control problem. The idea of studying a weak formulation stems from deterministic control problems, in which one needs to consider a family of optimization problems with different initial times and states. In the stochastic setting, in order to get deterministic initial condition on different initial time, we need to consider conditional probability spaces.
\end{remark}
For any fixed $s\in [0,T]$, let $\mathscr{P}_{[s,T]}$ denote the collection of all five-tuples stochastic systems
$\nu:=\left(\Omega,\mathscr{F},(\mathscr{F}_{t})_{t\in[s,T]},\mathbb{P},(W_{t})_{t\in[s,T]}\right)$, satisfying the following two conditions:
\begin{itemize}
\item $(\Omega,\mathscr{F},(\mathscr{F}_{t})_{t\in[s,T]},\mathbb{P})$ is a complete filtered probability space;
\item $(W_{t})_{t\in[s,T]}$ is a standard Brownian motion with respect to $(\mathscr{F}_{t})_{t\in[s,T]}$, defined on $(\Omega,\mathcal{F},\mathbb{P})$, with $W_{s}=0$, $\mathbb{P}$-a.$\,$s..
\end{itemize}
For any $\nu\in\mathscr{P}_{[s,T]}$, let
\begin{align*}
\mathcal{U}_{\nu}[s,T]:=\left\{U:[s,T]\times\Omega\rightarrow\mathscr{U}:\,\,U\text{ is }\,(\mathscr{F}_{t})_{t\in[s,T]}-\text{progressively measurable}\right\}.
\end{align*}
For $(s,x)\in [0,T]\times\overline{\mathscr{O}}$, $\nu\in\mathscr{P}_{[s,T]}$, consider the state equation \eqref{eq:StateEq} with initial condition
\begin{align}\label{eq:StateEqInits}
X_{s}=x=(p,z,\theta).
\end{align}
For any $U=(\Xi,L,C)\in\mathcal{U}_{\nu}[s,T]$, let
\begin{align}\label{eq:WeakCostFunt}
J_{\nu}(s,x,U)&:=\mathbb{E}^{s,x}\left(\int_{s}^{\tau_{\mathscr{O}}}\mathcal{L}(X_{t},U_{t})\,dt-P_{\tau_{\mathscr{O}}}Z_{\tau_{\mathscr{O}}}\right),\\
\label{eq:WeakValFunt1}
V_{\nu}(s,x)&:=\inf_{U\in\mathcal{U}_{\nu}[s,T]}J_{\nu}(s,x,U),\\
\label{eq:WeakValFunt2}
V(s,x)&:=\inf_{\nu\in\mathscr{P}_{[s,T]}}V_{\nu}(s,x),
\end{align}
where
\begin{align*}
\tau_{\mathscr{O}}=\tau_{\mathscr{O}}(s,x;U):=\inf\{t\geq s:\,\,X_{t}\not\in\mathscr{O},\,X_{s}=x\}\wedge T.
\end{align*}

\noindent
\textbf{Problem (SC').}\,\,\,Given any $(s,x)\in[0,T]\times\overline{\mathscr{O}}$, minimize (\ref{eq:WeakCostFunt}) over all $U\in\mathcal{U}_{\nu}[s,T]$ and all $\nu\in\mathscr{P}_{[s,T]}$. That is, find a five-tuples $\nu^{\ast}\in\mathscr{P}_{[s,T]}$ and $U^{\ast}\in\mathcal{U}_{\nu^{\ast}}[s,T]$, such that
\begin{align*}
J_{\nu}(s,x,U^{\ast})=V(s,x).
\end{align*}

Throughout this manuscript, various technical assumptions will be in order.
\begin{assumption}\label{assump:Coefficients}
The following standard assumptions are made on the coefficients of the state equations \eqref{eq:StateEqP}$-$\eqref{eq:StateEqTheta}.
\begin{itemize}
\item [(i)] There exists a constant $K>0$, such that for any $(\xi,\ell,c)\in\mathscr{U}$, $p,\tilde{p}\in[R,\infty)$ and any $t\in[0,T]$,
    \begin{align*}
    \left|b(t,p,\xi,\ell,c)-b(t,\tilde{p},\xi,\ell,c)\right|+\left|\Sigma(t,p,\xi,\ell,c)-\Sigma(t,\tilde{p},\xi,\ell,c)\right|&\leq K|p-\tilde{p}|,\\
    \left|b(t,p,\xi,\ell,c)\right|+\left|\Sigma(t,p,\xi,\ell,c)\right|&\leq K(1+p).
    \end{align*}
\item [(ii)] $b$ and $\Sigma$ are continuous on $[0,T]\times[R,\infty)\times\mathscr{U}$.
\item [(iii)] For any fixed $(\xi,\ell,c)\in\mathscr{U}$, $b(\cdot,\cdot,\xi,\ell,c),\,\sigma(\cdot,\cdot,\xi,\ell,c)\in C^{1,2}(\mathbb{R}^{+}\times[R,\infty))$.
\item [(iv)] $\vartheta,\,\sigma\in C^{1}([0,T])$.
\end{itemize}
\end{assumption}
\begin{assumption}\label{assump:SmallBallCond}
There exists a function $\psi:\mathbb{R}^{+}\rightarrow\mathbb{R}^{+}$, right-differentiable at the origin with $\psi(0)=\psi'_{+}(0)=0$ (where $\psi'_{+}(0)$ denotes the right derivative of $\psi$ at $0$), which is non-decreasing in a neighborhood of the origin and such that, for any $\varepsilon>0$, any $(s,x)\in[0,T]\times\overline{\mathscr{O}}$ with $\text{dist}(x,\partial\mathscr{O}\setminus\{z=0\})\leq\varepsilon$, $\nu\in\mathscr{P}_{[s,T]}$, $U\in\mathcal{U}_{\nu}[s,T]$,
\begin{align*}
\mathbb{E}^{s,x}\left(\tau_{\mathscr{O}}\right)\leq\psi(\varepsilon).
\end{align*}
\end{assumption}
\begin{remark}\label{rem:MainAssump}
Assumption \ref{assump:Coefficients}$-$(i) and (iv) ensure the existence of a unique strong solution to the SDE \eqref{eq:StateEqP}, while (ii)$-$(iv) of Assumption \ref{assump:Coefficients} and Assumption \ref{assump:SmallBallCond} are technical conditions for later proofs. Specifically, Assumption \ref{assump:Coefficients}$-$(ii), (iii), and (iv) are essential to prove the existence of smooth solutions to the HJB equation when uniform ellipticity is valid, and Assumption \ref{assump:SmallBallCond} is used to approximate the value function by smooth solutions to uniformly elliptic HJB equations.
\end{remark}
Before moving forward, we first verify that the expectation in \eqref{eq:WeakCostFunt} is finite, making the problem well defined.
\begin{lemma}\label{lem:WellDefStoContProb}
Under Assumption \ref{assump:Coefficients}$-$(i) and (iv), for any $(s,x)\in[0,T]\times\overline{\mathscr{O}}$, any $\nu\in\mathscr{P}_{[s,T]}$, and any $U\in\mathcal{U}_{\nu}[s,T]$, there exists a constant $\widetilde{K}>0$, which depends on $K$, $\kappa$, $\mathscr{L}$, $\mathscr{C}$, $\mathscr{N}$, $A$, $\alpha$, $\beta$, $\gamma$, $\varrho$ and $T$, such that
\begin{align*}
\left|J_{\nu}(s,x,U)\right|\leq\widetilde{K}\left(1+z+z^{2}+p^{2}\right),
\end{align*}
where $x=(p,z,\theta)$. In particular, \eqref{eq:WeakCostFunt}-\eqref{eq:WeakValFunt2} are well defined.
\end{lemma}
\noindent
\textbf{Proof:} Given any $(s,x)\in\overline{\mathscr{O}}_{T}$ with $x=(p,z,\theta)$, $\nu\in\mathscr{P}_{[s,T]}$ and $U\in\mathcal{U}_{\nu}[s,T]$, and since $(Z_{t})_{t\in[s,T]}$ is a $\mathbb{P}$-martingale with respect to $(\mathscr{F}_{t})_{t\in[s,T]}$,
\begin{align}\label{eq:EstCostFunt}
\left|J_{\nu}(s,x,U)\right|\leq\mathbb{E}^{s,x}\left(\int_{s}^{T}\kappa L_{t}^{\gamma}Z_{t}\,dt\right)+\mathbb{E}^{s,x}\left(P_{\tau_{\mathscr{O}}Z_{\tau_{\mathscr{O}}}}\right)\leq\kappa\mathscr{L}^{\gamma}zT+\mathbb{E}^{s,x}\left(P_{\tau_{\mathscr{O}}Z_{\tau_{\mathscr{O}}}}\right).
\end{align}
With the help of Assumption \ref{assump:Coefficients}$-$(i) as well as (D.5) in~\cite{FlemingSoner:2006},
\begin{align}\label{eq:EstPayoffTau}
\mathbb{E}^{s,x}\left(P^{2}_{\tau_{\mathscr{O}}}\right)\leq\mathbb{E}^{s,x}\left(\bigg(\sup_{t\in[s,T]}|P_{t}|\bigg)^{2}\right)\leq K_{1}\left(1+p^{2}\right)e^{K_{1}T},
\end{align}
where $K_{1}>0$ is a constant depending on $K$ and $T$. Moreover, since for any $t\in[s,T]$,
\begin{align*}
\mathbb{E}^{s,x}\left(Z^{2}_{t}\right)&\leq 2z^{2}+2\varrho^{-2}\int_{s}^{t}\mathbb{E}^{s,x}\left[Z_{r}^{2}\left(\Theta_{r}+A\,\Xi_{r}^{\alpha}L_{r}^{\beta}-C_{r}\right)^{2}\right]dr\\
&\leq 2z^{2}+2\left(H+1+A\mathscr{N}^{\alpha}\mathscr{L}^{\beta}+\mathscr{C}\right)^{2}\varrho^{-2}\int_{s}^{t}\mathbb{E}^{s,x}\left(Z^{2}_{r}\right)dr,
\end{align*}
and so by Gronwall's Inequality,
\begin{align*}
\mathbb{E}^{s,x}\left(Z^{2}_{t}\right)&\leq 2z^{2}\left(1+\exp\left(2\left(H+1+A\mathscr{N}^{\alpha}\mathscr{L}^{\beta}+\mathscr{C}\right)^{2}\varrho^{-2}T\right)\right).
\end{align*}
Therefore,
\begin{align}
\mathbb{E}^{s,x}\!\left(Z^{2}_{\tau_{\mathscr{O}}}\right)&\leq 2z^{2}+2\varrho^{-2}\int_{s}^{T}\mathbb{E}^{s,x}\left[Z^{2}_{r}\left(\Theta_{r}+A\,\Xi_{r}^{\alpha}L_{r}^{\beta}-C_{r}\right)^{2}\right]dr\nonumber\\
&\leq 2z^{2}+2\varrho^{-2}\left(H+1+A\mathscr{N}^{\alpha}\!\mathscr{L}^{\beta}+\mathscr{C}\right)^{2}\int_{s}^{T}\mathbb{E}^{s,x}\left(Z^{2}_{r}\right)dr\nonumber\\
\label{eq:EstGirDenTau} &\leq 2z^{2}+\frac{4z^{2}T}{\varrho^{2}}\!\left(H\!+\!1\!+\!A\mathscr{N}^{\alpha}\!\mathscr{L}^{\beta}\!+\!\mathscr{C}\right)^{2}\!\left(\!1\!+\!\exp\left(\frac{2T}{\varrho^{2}}\!\left(H\!+\!1\!+\!A\mathscr{N}^{\alpha}\!\mathscr{L}^{\beta}\!+\!\mathscr{C}\right)^{2}\right)\!\right).
\end{align}
Combining \eqref{eq:EstCostFunt}, \eqref{eq:EstPayoffTau}, and \eqref{eq:EstGirDenTau} completes the proof.\hfill $\Box$

\subsection{The HJB Equation and Viscosity Solutions}\label{subsec:HJBVisSol}

Let $\mathscr{S}^{n}$ be the set of all $n\times n$ symmetric matrices, equipped with its usual order. That is, for $G_{1},G_{2}\in\mathscr{S}^{n}$, $G_{1}\leq G_{2}$ if and only if $G_{2}-G_{1}$ is non-negative definite. Let $\mathscr{S}_{+}^{n}\subseteq\mathscr{S}^{n}$ be the subset of nonnegative-definite $n\times n$ matrices. For $s\in [0,T]$, $x=(p,z,\theta)\in\overline{\mathscr{O}}$ and $u=(\xi,\ell,c)\in\mathscr{U}$, $M\in\mathscr{S}_{+}^{3}$ and $y\in\mathbb{R}^{3}$, define the Hamiltonian
\begin{align*}
\mathcal{H}(s,x,y,M)=\sup_{u\in\mathscr{U}}\left(-\vec{f}(s,x,u)\cdot y-\frac{1}{2}\text{tr}\left(\vec{a}(s,x,u)M\right)-\mathcal{L}(x,u)\right).
\end{align*}
The HJB equation associated with the stochastic control problem \eqref{eq:WeakCostFunt}-\eqref{eq:WeakValFunt2} is
\begin{align}\label{eq:HJB}
-\frac{\partial V}{\partial s}+\mathcal{H}(s,x,D_{x}V,D_{x}^{2}V)=0,\quad (s,x)\in Q_{T}:=[0,T)\times\mathscr{O},
\end{align}
with terminal/boundary condition
\begin{align}\label{eq:HJBTermBound}
V(s,x)=-pz,\quad (s,x)\in\partial^{*}Q_{T}:=\left([0,T]\times\partial\mathscr{O}\right)\cup\left(\{T\}\times\mathscr{O}\right).
\end{align}

By standard stochastic control theory, the value function \eqref{eq:WeakValFunt2} is expected to be a classical solution of the HJB equation \eqref{eq:HJB} with terminal/boundary condition \eqref{eq:HJBTermBound}, provided that the following \emph{uniformly elliptic} condition holds: there exists a constant $\lambda_{0}>0$, such that
\begin{align}\label{eq:UnifEllip}
\sum_{i,j=1}^{3}a_{ij}(s,x,u)w_{i}w_{j}\geq\lambda_{0}|w|^{2},\quad\text{for all }\,(s,x)\in\overline{Q}_{T},\,\,\,u\in\mathscr{U},\,\,\,w\in\mathbb{R}^{3}.
\end{align}
Unfortunately, our stochastic control system does not satisfy \eqref{eq:UnifEllip}. In particular, the matrix $\vec{a}(t,x,u)$ is not even positive definite. Hence, we can only connect our value function with the HJB equation via a viscosity solution. Throughout, let $C(\overline{Q}_{T})$ be the set of continuous functions on $\overline{Q}_{T}$, and let $C^{1,2}(\overline{Q}_{T})$ be the set of all functions $f$ whose partial derivatives $(\partial f/\partial t)$, $(\partial f/\partial x_{i})$, $(\partial f/\partial x_{i}\partial x_{j})$ exist and are continuous on $\overline{Q}_{T}$. Now, recall (cf.~\cite[Definition II.4.1]{FlemingSoner:2006}):
\begin{definition}\label{def:VisSol}
A function $v\in C(\overline{Q}_{T})$ is called a \emph{viscosity subsolution} to \eqref{eq:HJB} with terminal condition \eqref{eq:HJBTermBound} if,
\begin{align}\label{eq:VisSubBound}
v(s,x)\leq -pz,\quad\text{for any }\,(s,x)\in\partial^{*}Q_{T},
\end{align}
and if, for any $\varphi\in C^{1,2}(\overline{Q}_{T})$ such that $v-\varphi$ attains a local maximum at some $(\bar{s},\bar{x})\in Q_{T}$,
\begin{align*}
-\varphi_{s}(\bar{s},\bar{x})+\mathcal{H}\left(\bar{s},\bar{x},D_{x}\varphi(\bar{s},\bar{x}),D_{xx}\varphi(\bar{s},\bar{x})\right)\leq 0.
\end{align*}
Similarly, a function $v\in C(\overline{Q}_{T})$ is called a \emph{viscosity supersolution} of \eqref{eq:HJB} with terminal condition \eqref{eq:HJBTermBound} if,
\begin{align}\label{eq:VisSuperBound}
v(s,x)\geq -pz,\quad\text{for any }\,(s,x)\in\partial^{*}Q_{T},
\end{align}
and if, for any $\varphi\in C^{1,2}(\overline{Q}_{T})$ such that $v-\varphi$ attains a local minimum at some $(\bar{s},\bar{x})\in Q_{T}$,
\begin{align*}
-\varphi_{s}(\bar{s},\bar{x})+\mathcal{H}\left(\bar{s},\bar{x},D_{x}\varphi(\bar{s},\bar{x}),D_{xx}\varphi(\bar{s},\bar{x})\right)\geq 0.
\end{align*}
$v$ is called a \emph{viscosity solution} to \eqref{eq:HJB} with terminal condition \eqref{eq:HJBTermBound} if it is both a
viscosity subsolution and a viscosity supersolution.
\end{definition}

In studying viscosity solutions of a second-order parabolic HJB equation, an equivalent definition in terms of second-order sub-differentials and super-differentials is useful (cf.~\cite[Section 8]{CrandallIshiiLions:1992} and~\cite[Definition V.4.1 \& V.4.2]{FlemingSoner:2006}).
\begin{definition}\label{def:SubSuperDiff}
Let $v\in C(\overline{Q}_{T})$.
\begin{itemize}
\item [(i)] The set of (parabolic) second-order \emph{super-differentials} of $v$ at $(s,x)\in Q_{T}$ is
\begin{align*}
\mathcal{D}^{(1,2)}_{+}\!v(s,x)\!:=\!\left\{\!(q,r,G)\!\in\!\mathbb{R}\!\times\!\mathbb{R}^{3}\!\times\!\mathscr{S}^{3}\!:\!v(s\!+\!h,x\!+\!y)\!-\!v(s,x)\!\leq\!qh\!+\!r\!\cdot\!y\!+\!\frac{1}{2}y\!\cdot\!Gy\!+\!o(|h|\!+\!|y|^{2})\!\right\}.
\end{align*}
\item [(ii)] The set of (parabolic) second-order \emph{sub-differentials} of $v$ at $(t,x)\in Q_{T}$ is
\begin{align*}
\mathcal{D}^{(1,2)}_{-}\!v(s,x)\!&:=\!-D^{(1,2)}_{+}(-v)(s,x)\\
&\,=\left\{\!(q,r,G)\!\in\!\mathbb{R}\!\times\!\mathbb{R}^{3}\!\times\!\mathscr{S}^{3}\!:\!v(s\!+\!h,x\!+\!y)\!-\!v(s,x)\!\geq\!qh\!+\!r\!\cdot\!y\!+\!\frac{1}{2}y\!\cdot\!Gy\!+\!o(|h|\!+\!|y|^{2})\!\right\}.
\end{align*}
\item [(iii)] The closure of the set of sub- and super-differentials of $v$ at $(t,x)\in Q_{T}$ are
\begin{align*}
\overline{\mathcal{D}}^{(1,2)}_{\pm}v(s,x)&:=\left\{(q,r,G)\in\mathbb{R}\times\mathbb{R}^{3}\times\mathscr{S}^{3}:\,\exists(s_{n},x_{n})\in Q_{T},\,(q_{n},r_{n},G_{n})\in D^{(1,2)}_{\pm}v(s_{n},x_{n}),\right.\\
&\qquad\qquad\qquad\qquad\qquad\qquad\,\,\,\left.\lim_{n\rightarrow\infty}(s_{n},x_{n})=(s,x),\,\lim_{n\rightarrow\infty}(q_{n},r_{n},G_{n})=(q,r,G)\right\}.
\end{align*}
\end{itemize}
\end{definition}
It follows from Definition \ref{def:SubSuperDiff} that, if $\varphi\in C^{1,2}(Q_{T})$, then
\begin{align*}
\mathcal{D}^{(1,2)}_{\pm}(v-\varphi)(s,x)=\left\{\left(q-\varphi_{s}(s,x),r-D_{x}\varphi(s,x),G-D_{xx}\varphi(s,x)\right):(q,r,G)\in\mathcal{D}^{(1,2)}_{\pm}v(s,x))\right\}.
\end{align*}
The same statement holds for $\overline{\mathcal{D}}^{(1,2)}_{\pm}$. Moreover, the characterizations of the
second-order sub and superdifferentials in Definition \ref{def:SubSuperDiff} yield
\begin{align}\label{eq:SuperDiff}
-q+\mathcal{H}(s,x,r,G)&\leq 0,\quad\text{for any }\,(q,r,G)\in\overline{\mathcal{D}}^{(1,2)}_{+}v(s,x),\\
\label{eq:Subdiff} -q+\mathcal{H}(s,x,r,G)&\geq 0,\quad\text{for any }\,(q,r,G)\in\overline{\mathcal{D}}^{(1,2)}_{-}v(s,x).
\end{align}
The above inequalities form an equivalent requirement for viscosity sub- and super-solutions. Towards obtaining this equivalence, we start by stating the following lemma whose proof can be found in \cite[Lemma V.4.1]{FlemingSoner:2006}.
\begin{lemma}\label{lem:RepSubSuperDiff}
Let $v\in C(\overline{Q}_{T})$, and let $(s,x)\in Q_{T}$. Then, $(q,r,G)\in\mathcal{D}^{(1,2)}_{+}v(s,x)$ if and only if there exists $\tilde{v}\in C^{1,2}(\overline{Q}_{T})$, such that
\begin{align}\label{eq:RepSubSuperDiff}
\left(\tilde{v}_{s}(s,x),D_{x}\tilde{v}(s,x),D_{xx}\tilde{v}(s,x)\right)=(q,r,G),
\end{align}
such that $v-\tilde{v}$ attains its maximum at $(s,x)$. Similarly, $(q,r,G)\in\mathcal{D}^{(1,2)}_{-}v(s,x)$ if and only if there exists $\hat{v}\in C^{1,2}(\overline{Q}_{T})$ satisfying \eqref{eq:RepSubSuperDiff}, such that $v-\hat{v}$ attains its minimum at $(s,x)$.
\end{lemma}
An immediate corollary to the above result is the following equivalent definition of viscosity solution for the (second-order) HJB equation \eqref{eq:HJB} with boundary/terminal condition \eqref{eq:HJBTermBound} (cf.~\cite[Proposition 4.1]{FlemingSoner:2006}).
\begin{proposition}\label{prop:EquCharVisSols}
$v\in C(\overline{Q}_{T})$ is a viscosity subsolution of \eqref{eq:HJB} with terminal/boundary condition \eqref{eq:HJBTermBound} if and only if \eqref{eq:SuperDiff} holds for all $(t,x)\in Q_{T}$, and (\ref{eq:VisSubBound}) holds for all $(t,x)\in \partial^{*}Q_{T}$. Similarly, $v\in C(\overline{Q}_{T})$ is a viscosity supersolution of \eqref{eq:HJB} with terminal/boundary condition \eqref{eq:HJBTermBound} if and only if \eqref{eq:SuperDiff} holds for all $(t,x)\in Q_{T}$, and \eqref{eq:VisSuperBound} holds for all $(t,x)\in\partial^{*}Q_{T}$.
\end{proposition}

\section{Existence of Viscosity Solution}\label{sec:Existence}

The main goal of the section is to verify that the value function $V$, as given in \eqref{eq:WeakValFunt2}, is indeed a viscosity solution of the HJB equation \eqref{eq:HJB} with terminal/boundary condition \eqref{eq:HJBTermBound}. The proof will proceed in three steps. In Section \ref{subsec:ContValFunt}, we first justify the joint continuity of the value function. In Section \ref{subsec:DynProgPrin}, we investigate the so-called Dynamic Programming Principle for the value function. The main difficulty in both sections stems from that, for different time values $t$, the value function $V$ is defined based on different probability spaces. Therefore, one cannot prove the joint continuity or the Dynamic Programming Principle by direct estimations of expectations. Finally, in Section \ref{subsec:ExistVisSols}, based on the joint continuity and the Dynamic Programming Principle, the value function is shown to satisfy Definition \ref{def:VisSol}.

\subsection{Continuity of the Value Function}\label{subsec:ContValFunt}

The main tool in verifying the joint continuity of the value function \eqref{eq:WeakValFunt2} is a perturbation method similar to~\cite[Lemma IV.7.1 \& Theorem IV.7.2]{FlemingSoner:2006}. More precisely, we will approximate the HJB equation \eqref{eq:HJB} by adding a small perturbation so that the uniform ellipticity \eqref{eq:UnifEllip} is satisfied. Moreover, we will restrict the domain of the state equation to a compact subspace on which the perturbed HJB equation has a unique classical solution. The value function \eqref{eq:WeakValFunt2} will then be identified as the \emph{uniform} limit, over all possible controls, of this classical solution by taking the perturbation to zero and the bounded domain to the original half-unbounded domain. Note, however, that our stochastic control problem lies on a half-unbounded domain $\mathscr{O}$, rather than the whole Euclidean space as in~\cite[Section IV.7]{FlemingSoner:2006}, and thus the exit time $\tau$ is considered in \eqref{eq:WeakCostFunt} instead of the terminal time $T$, which greatly increases the difficulty in the perturbation method.

Before stating the main theorem of this section, we first introduce some notations and one extra technical assumption. For any $\rho>R$, let $(R,\rho)\times\left(\rho^{-1},\rho\right)\times\left(-H+\rho^{-1},H-\rho^{-1}\right)\subset\mathscr{O}_{\rho}\subset\mathscr{O}$, such that $\partial\mathscr{O}_{\rho}\in C^{3}(\mathbb{R}^{3})$ and that $\mathscr{O}_{\rho}$ is increasing in $\rho$. Then as $\rho\rightarrow\infty$, $\mathscr{O}_{\rho}\uparrow\mathscr{O}$. Let $\alpha_{\rho}:\overline{\mathscr{O}}\rightarrow[0,1]$ be such that $\alpha_{\rho}\in C^{\infty}(\overline{\mathscr{O}})$, $\alpha_{\rho}(x)=1$ for $x\in\overline{\mathscr{O}}_{\rho}$, and that $\alpha_{\rho}(x)=0$ for
$x\in\overline{\mathscr{O}}\setminus\overline{\mathscr{O}}_{\rho+1}$. Also, for $s\in[0,T]$, $x=(p,z,\theta)\in\overline{\mathscr{O}}$ and $u=(\xi,\ell,c)\in\mathscr{U}$, denote
\begin{align*}
\vec{f}_{\rho}(s,x,u):=\vec{f}(s,x,u)\alpha_{\rho}(x),\quad\vec{\sigma}_{\rho}(s,x,u):=\vec{\sigma}(s,x,u)\alpha_{\rho}(x),\quad\mathcal{L}_{\rho}(x,u):=\mathcal{L}(x,u)\alpha_{\rho}(x).
\end{align*}
Next, for any fixed $s\in[0,T]$, let $\widetilde{\mathscr{P}}_{[s,T]}$ denote the collection of all six-tuple stochastic systems $\mu:=(\Omega,\mathscr{F},(\mathscr{F}_{t})_{t\in[s,T]},\mathbb{P},W,\widetilde{W})$, where
\begin{itemize}
\item $(\Omega,\mathscr{F},(\mathscr{F}_{t})_{t\in[s,T]},\mathbb{P})$ is a complete filtered probability space;
\item $W:=(W_{t})_{t\in[s,T]}$ is a one-dimensional standard Brownian motion with $W_{s}=0$, $\mathbb{P}$-a.$\,$s.;
\item $\widetilde{W}:=(\widetilde{W}_{t})_{t\in[s,T]}$ is a three-dimensional standard Brownian motion, independent of $W$, with $\widetilde{W}_{s}=0$, $\mathbb{P}$-a.$\,$s..
\end{itemize}
Under each $\mu\in\widetilde{\mathscr{P}}_{[s,T]}$, let $\mathcal{V}_{\mu}[s,T]$ be the collection of progressively measurable processes on $\Omega\times[s,T]$, taking values in $\mathscr{U}$.

Fix $\epsilon>0$, for any $\mu\in\widetilde{\mathscr{P}}_{[s,T]}$, any $U=(U_{t})_{t\in[s,T]}\in\mathcal{V}_{\mu}[s,T]$, and any $n\in\mathbb{N}$, consider the state equation for $X^{(\rho,n)}=(X_{t}^{(\rho,n)})_{t\in[s,T]}$, where $X_{t}^{(\rho,n)}:=(P_{t}^{(\rho,n)},Z_{t}^{(\rho,n)},\Theta_{t}^{(\rho,n)})$, $t\in[s,T]$,
\begin{align}\label{eq:StateEqRhon}
dX_{t}^{(\rho,n)}=\vec{f}_{\rho}\left(t,X_{t}^{(\rho,n)},U_{t}\right)dt+\vec{\sigma}_{\rho}\left(t,X_{t}^{(\rho,n)},U_{t}\right)dW_{t}+\epsilon^{n}I_{3}\,d\widetilde{W}_{t},
\end{align}
with initial condition $X_{s}^{(\rho,n)}=x$, where $I_{3}$ denote the $3\times 3$ identity matrix. Also, let
\begin{align}\label{eq:RhonCostFunt}
J_{\rho,\mu}^{(n)}(s,x;U)&:=\mathbb{E}^{s,x}\left(\int_{s}^{\tau_{\rho,n}}\mathcal{L}_{\rho}\left(X_{t}^{(\rho,n)},U_{t}\right)dt-P_{\tau_{\rho,n}}^{(\rho,n)}Z_{\tau_{\rho,n}}^{(\rho,n)}\right),\\
\label{eq:RhonValFunt1} V_{\rho,\mu}^{(n)}(s,x)&:=\inf_{U\in\mathcal{V}_{\mu}[s,T]}J_{\rho,\mu}^{(n)}(s,x;U),\\
\label{eq:RhonValFunt2} V_{\rho}^{(n)}(s,x)&:=\inf_{\mu\in\widetilde{\mathscr{P}}_{[s,T]}}V_{\rho,\mu}^{(n)}(s,x),
\end{align}
where
\begin{align*}
\tau_{\rho,n}:=\tau_{\rho,n}(s,x)=\inf\left\{t\geq s:\,X_{t}^{(\rho,n)}\notin\mathscr{O}_{\rho},\,X_{s}^{(\rho,n)}=x\right\}\wedge T.
\end{align*}
\begin{remark}\label{rem:EquiStoSym}
Above, $V_{\rho,\mu}^{(n)}(s,y)$ is defined on a six-tuple stochastic system $\mu$ instead of on a five-tuple stochastic system $\nu$, where $V_{\nu}$ (see \eqref{eq:WeakValFunt1}) is defined. However, given $\mu=(\Omega,\mathscr{F},(\mathscr{F}_{t})_{t\in[s,T]},\mathbb{P},W,\widetilde{W})$, by setting $\nu=(\Omega,\mathscr{F},(\mathscr{F}_{t})_{t\in[s,T]},\mathbb{P},W)$, any $U\in\mathcal{V}_{\mu}[s,T]$ is also in $\mathcal{U}_{\nu}[s,T]$. On the other hand, given a five-tuple $\nu\in\mathscr{P}_{[s,T]}$, consider another three five-tuples $\nu_{i}=(\Omega_{i},\mathscr{F}^{(i)},(\mathcal{F}_{t}^{(i)})_{t\in[s,T]},\mathbb{P}_{i},W^{(i)})$, $i=1,2,3$. For $(\omega,\omega_{1},\omega_{2},\omega_{3})\in\widehat{\Omega}:=\Omega\times\Omega_{1}\times\Omega_{2}\times\Omega_{3}$, set
\begin{align*}
\widehat{W}_{t}(\omega,\omega_{1},\omega_{2},\omega_{3})=W_{t}(\omega),\quad\widetilde{W}_{t}(\omega,\omega_{1},\omega_{2},\omega_{3})=\left(W_{t}^{(1)}(\omega_{1}),W_{t}^{(2)}(\omega_{2}),W_{t}^{(3)}(\omega_{3})\right)^{T}.
\end{align*}
Hence,
\begin{align*}
\hat{\mu}:=\left(\widehat{\Omega},\mathscr{F}\otimes\left(\bigotimes_{i=1}^{3}\mathscr{F}^{(i)}\right),\left(\mathscr{F}_{t}\otimes\left(\bigotimes_{i=1}^{3}\mathscr{F}_{t}^{(i)}\right)\right)_{t\in[s,T]},\mathbb{P}\otimes\left(\bigotimes_{i=1}^{3}\mathbb{P}_{i}\right),\widehat{W},\widetilde{W}\right),
\end{align*}
is an element in $\widetilde{\mathscr{P}}_{[s,T]}$. Thus, any $U\in\mathcal{U}_{\nu}[s,T]$ can also be regarded as an element in $\mathcal{V}_{\hat{\mu}}[s,T]$.
\end{remark}
To proceed to the proof of the main theorem, we need the following technical assumption, which is the analog of Assumption \ref{assump:SmallBallCond} for $X^{(\rho,n)}$.
\begin{assumption}\label{assump:SmallBallCondRhon}
Let $\psi$ be as in Assumption \ref{assump:SmallBallCond}. For any $\rho>R$, $\varepsilon>0$, any $(s,x)\in[0,T]\times\overline{\mathscr{O}}_{\rho}$ with $\text{dist}(x,\partial\mathscr{O}_{\rho})\leq\varepsilon$, any $\mu\in\widetilde{\mathscr{P}}_{[s,T]}$, $U\in\mathcal{V}_{\mu}[s,T]$, and $n\in\mathbb{N}$,
\begin{align*}
\mathbb{E}^{s,x}\left(\tau_{\rho,n}\right)\leq\psi(\varepsilon).
\end{align*}
\end{assumption}
Next, the HJB equation associated with \eqref{eq:RhonCostFunt}-\eqref{eq:RhonValFunt2} is
\begin{align}\label{eq:RhonHJB}
-\frac{\partial V_{\rho}^{(n)}}{\partial
s}+\mathcal{H}_{\rho}^{(n)}\left(s,x,D_{x}V_{\rho}^{(n)},D_{x}^{2}V_{\rho}^{(n)}\right)=0,\quad (s,x)\in Q_{\rho,T}:=[0,T)\times\mathscr{O}_{\rho},
\end{align}
with terminal/boundary condition
\begin{align}\label{eq:RhonHJBTermBound}
V_{\rho}^{(n)}(s,x)=-pz,\quad (s,x)\in\partial^{*}Q_{\rho,T},
\end{align}
where, for $M\in\mathscr{S}_{+}^{3}$ and $y\in\mathbb{R}^{3}$,
\begin{align*}
\mathcal{H}_{\rho}^{(n)}(s,x,y,M)&:=\sup_{u\in\mathscr{U}}\left(-\vec{f_{\rho}}(t,x,u)\cdot y-\frac{1}{2}\text{tr}(\vec{a}_{\rho}^{(n)}(s,x,u)M)-L_{\rho}(x,u)\right),
\end{align*}
and where $\vec{a}_{\rho}^{(n)}(s,x,u):=2(\vec{\sigma}_{\rho}(s,x,u)+\epsilon^{n}I_{3})(\vec{\sigma}_{\rho}(s,x,u)+\epsilon^{n}I_{3})^{T}$. Note that we treat $\vec{\sigma}_{\rho}$ as a $3\times 3$ matrix with the last two columns both identically zero, when taking the summation with $\epsilon^{n}I_{3}$.
\begin{theorem}\label{thm:JointCont}
Under Assumptions \ref{assump:Coefficients}, Assumption \ref{assump:SmallBallCond}, and Assumption \ref{assump:SmallBallCondRhon}, the value function $V$, as given in \eqref{eq:WeakValFunt2}, is continuous on $\overline{Q}_{T}$. Moreover, for every $(s,y)\in\overline{Q}_{T}$, $V(s,y)=V_{\nu}(s,y)$, for all $\nu\in\mathscr{P}_{[s,T]}$.
\end{theorem}

\noindent
\textbf{Proof: Step 1.} We first consider the stochastic control problem \eqref{eq:RhonCostFunt}-\eqref{eq:RhonValFunt2}. Note that the SDE \eqref{eq:StateEqRhon} satisfies the uniform ellipticity condition \eqref{eq:UnifEllip}. By standard stochastic control theory (cf.~\cite[Theorem IV.4.1]{FlemingSoner:2006}, the conditions therein are satisfied from Assumptions \ref{assump:Coefficients}) , the HJB equation \eqref{eq:RhonHJB} has a unique solution $\mathcal{W}^{(n)}_{\rho}\in C^{1,2}(Q_{\rho,T})\cap C(\overline{Q}_{\rho,T})$ with terminal/boundary condition \eqref{eq:RhonHJBTermBound}. We will verify the joint continuity of $V_{\rho}^{(n)}$ by showing that $\mathcal{W}^{(n)}_{\rho}(s,x)=V_{\rho,\mu}^{(n)}(s,x)$ for any $\mu\in\widetilde{\mathscr{P}}_{[s,T]}$ and  $(s,x)\in\overline{Q}_{\rho,T}$, which also implies that $V_{\rho,\mu}^{(n)}=V_{\rho}^{(n)}$ for all $\mu\in\widetilde{\mathscr{P}}_{[s,T]}$. In the following proof, we fix $s\in[0,T)$ and $x\in\overline{\mathscr{O}}\setminus\{z=0\}$, and choose $\rho>R$ large enough so that $x\in\overline{\mathscr{O}}_{\rho}$.

To start with, choose $\delta>0$ so that $\rho-\delta>R$ and $x\in\overline{\mathscr{O}}_{\rho-\delta}$, then $\partial\mathcal{W}^{(n)}_{\rho}/\partial s$, $D_{x}\mathcal{W}^{(n)}_{\rho}$ and $D_{xx}\mathcal{W}^{(n)}_{\rho}$ are all uniformly continuous on $\overline{Q}_{\rho-\delta,T}:=[0,T-\delta]\times\overline{\mathscr{O}}_{\rho-\delta}$. Hence, for any $\varepsilon>0$, there exists $\kappa_{1}>0$, such that for any $(\tilde{s},\tilde{x}),(\hat{s},\hat{x})\in\overline{Q}_{\rho-\delta,T}$ with $|\tilde{s}-\hat{s}|<\kappa_{1}$ and $|\tilde{x}-\hat{x}|<\kappa_{1}$,
\begin{align}\label{eq:UnifContWRhon}
\left|\mathscr{A}_{\rho,u}^{(n)}\mathcal{W}^{(n)}_{\rho}(\tilde{s},\tilde{x})+\mathcal{L}_{\rho}(\tilde{x},u)-\mathscr{A}_{\rho,u}^{(n)}\mathcal{W}^{(n)}_{\rho}(\hat{s},\hat{x})-\mathcal{L}_{\rho}(\hat{x},u)\right|<\frac{\varepsilon}{4T},\quad\text{for all }\,u\in\mathscr{U},
\end{align}
where $\mathscr{A}_{\rho,u}^{(n)}$ is a parabolic operator defined by:
\begin{align*}
\mathscr{A}_{\rho,u}^{(n)}\mathcal{W}^{(n)}_{\rho}:=\frac{\partial\mathcal{W}_{\rho}^{(n)}}{\partial s}+\frac{1}{2}\text{tr}\left(\vec{a}_{\rho}^{(n)}D_{xx}\mathcal{W}_{\rho}^{(n)}\right)+\vec{f}_{\rho}\cdot D_{x}\mathcal{W}_{\rho}^{(n)}.
\end{align*}
The HJB equation \eqref{eq:RhonHJB} can then be written as
\begin{align}\label{eq:RhonHJB2}
\min_{u\in\mathscr{U}}\left(\mathscr{A}_{\rho,u}^{(n)}\mathcal{W}^{(n)}_{\rho}(s,x,u)+\mathcal{L}_{\rho}(x,u)\right)=0.
\end{align}
Choose $M_{1}>0$ large so that $(T-s)/M_{1}<\min(\kappa_{1},1)$, and divide $[s,T-\delta)$ into $M_{1}$ subintervals $\mathcal{I}_{i}=[s_{i},s_{i+1})$, $i=1,\ldots,M_{1}$. Also, choose $M_{2}>0$ large and partition $\overline{\mathscr{O}}_{\rho-\delta}$ into disjoint Borel sets: $\overline{\mathscr{O}}_{\rho-\delta}=\mathcal{B}_{1}\cup\cdots\cup\mathcal{B}_{M_{2}}$, so that each $\mathcal{B}_{j}$ has the diameter no more than $\kappa_{1}/2$. Pick $x_{j}\in\mathcal{B}_{j}$. For each $i=1,\ldots,M_{1}$, $j=1,\ldots,M_{2}$, by \eqref{eq:RhonHJB2}, there exists $u_{ij}\in\mathscr{U}$, such that
\begin{align*}
\mathscr{A}_{\rho,u_{ij}}^{(n)}\mathcal{W}^{(n)}_{\rho}(s_{i},x_{j})+\mathcal{L}_{\rho}(x_{j},u_{ij})<\frac{\varepsilon}{4T},
\end{align*}
which, together with \eqref{eq:UnifContWRhon}, implies that, for $t\in I_{i}$, $|y-x_{j}|<\kappa_{1}$,
\begin{align}\label{eq:EstALRhonuij}
\mathscr{A}_{\rho,u_{ij}}^{(n)}\mathcal{W}^{(n)}_{\rho}(t,y)+\mathcal{L}_{\rho}(y,u_{ij})<\frac{\varepsilon}{2T}.
\end{align}

Pick an arbitrary $u_{0}\in\mathscr{U}$, and define the discrete Markov control policy $\underline{u}:=(\underline{u}_{1},\ldots,\underline{u}_{M})$ by
\begin{align*}
\underline{u}_{i}(y):=\left\{\begin{array}{ll} u_{ij} &\text{if }\,y\in\mathcal{B}_{j}\,\,\,\text{for some }\,j=1,\ldots,M_{2}, \\ u_{0} &\text{otherwise}. \end{array}\right.
\end{align*}
Define $\widehat{U}\in\mathcal{V}_{\mu}[s,T]$ and the solution $X^{(\rho,n)}$ to \eqref{eq:StateEqRhon} with $X^{(\rho,n)}_{s}=x$ and control $\widehat{U}$ such that
\begin{align*}
\widehat{U}_{t}=\underline{u}_{i}\left(X^{(\rho,n)}_{s_{i}}\right)\quad\text{for }\,t\in\mathcal{I}_{i},\quad i=1,\ldots,M_{1}.
\end{align*}
This can be done by induction on $i$, since for $t\in\mathcal{I}_{i}$, $X^{(\rho,n)}_{t}$ is the solution to \eqref{eq:StateEqRhon} with initial data $X^{(\rho,n)}_{s_{i}}$, and for $t\in[T-\delta,T]$, $X^{(\rho,n)}_{t}$ is the solution to \eqref{eq:StateEqRhon} with initial data $X^{(\rho,n)}_{t_{M+1}}$. In particular, $\widehat{U}_{t}=u_{ij}$ if $t\in\mathcal{I}_{i}$ and $X^{(\rho,n)}_{t}\in\mathcal{B}_{j}$.

By Dynkin's formula, for any $(\mathscr{F}_{t})_{t\in[s,T]}$-stopping time $\tau$,
\begin{align}\label{eq:DynkinWRhon1}
\mathcal{W}^{(n)}_{\rho}(s,x)&=\mathbb{E}^{s,x}\left(-\int_{s}^{\tau\wedge\tau_{\rho,n}^{(\delta)}}\mathscr{A}_{\rho,\widehat{U}}^{(n)}\,\mathcal{W}^{(n)}_{\rho}\!\left(t,X_{t}^{(\rho,n)}\right)dt+\mathcal{W}^{(n)}_{\rho}\!\left(\tau\wedge\tau_{\rho,n}^{(\delta)},X_{\tau\wedge\tau_{\rho,n}^{(\delta)}}^{(\rho,n)}\right)\right)\\
&=\mathbb{E}^{s,x}\left(\int_{s}^{\tau\wedge\tau_{\rho,n}^{(\delta)}}\mathcal{L}_{\rho}\!\left(X_{(\rho,n)}^{t},\widehat{U}_{t}\right)dt+\mathcal{W}^{(n)}_{\rho}\!\left(\tau\wedge\tau_{\rho,n}^{(\delta)},X_{\tau\wedge\tau_{\rho,n}^{(\delta)}}^{(\rho,n)}\right)\right)\nonumber\\
\label{eq:DynkinWRhon2} &\quad\,-\mathbb{E}^{s,x}\left(\int_{s}^{\tau\wedge\tau_{\rho,n}^{(\delta)}}\left(\mathscr{A}_{\rho,\widehat{U}}^{(n)}\,\mathcal{W}^{(n)}_{\rho}\!\left(t,X_{t}^{(\rho,n)}\right)+\mathcal{L}_{\rho}\!\left(X_{(\rho,n)}^{t},\widehat{U}_{t}\right)\right)dt\right).
\end{align}
where $\tau_{\rho,n}^{(\delta)}:=\inf\{t\geq s:\,X_{t}^{(\rho,n)}\notin\mathscr{O}_{\rho-\delta}\}\wedge(T-\delta)$. We need to estimate the second term in \eqref{eq:DynkinWRhon2}. To that effect, define
\begin{align*}
\Gamma=\left\{\omega\in\Omega:\,X_{t}^{(\rho,n)}(\omega)\in\overline{\mathscr{O}}_{\rho-\delta},\,\left|X_{t}^{(\rho,n)}(\omega)-X_{s_{i}}^{(\rho,n)}(\omega)\right|<\frac{\kappa_{1}}{2},\,t\in\mathcal{I}_{i},\,i=1,\ldots,M_{1}\right\}.
\end{align*}
By \eqref{eq:EstALRhonuij}, for any $t\in [s,T-\delta)$,
\begin{align}\label{eq:EstALRhonGamma}
\mathscr{A}_{\rho,\widehat{U}(\omega)}^{(n)}\,\mathcal{W}_{\rho}^{(n)}\!\left(t,X_{t}^{(\rho,n)}(\omega)\right)+\mathcal{L}_{\rho}\left(X_{t}^{(\rho,n)}(\omega),\widehat{U}_{t}(\omega)\right)<\frac{\varepsilon}{2T},\quad\omega\in\Gamma,
\end{align}
By Assumption \ref{assump:Coefficients}$-$(ii) and the very definition of $\vec{f}_{\rho}$ and $\vec{\sigma}_{\rho}$, the drift and volatility vector are both bounded. Hence, by (D.12) in~\cite{FlemingSoner:2006}, there exists a constant $D_{1}>0$, depending only on $\rho$, $T$ and $K$ (Assumption \ref{assump:Coefficients}$-$(i)), such that
\begin{align*}
\mathbb{P}^{s,x}\left(\max_{t\in\mathcal{I}_{i}}\left|X_{t}^{(\rho,n)}-X_{s_{i}}^{(\rho,n)}\right|\geq\frac{\kappa_{1}}{2}\right)\leq\kappa_{1}^{-4}D_{1}(s_{i+1}-s_{i})^{2},\quad\text{for all }\,i=1,\ldots,M_{1},
\end{align*}
from which it follows that
\begin{align}\label{eq:EstGammaComp}
\mathbb{P}^{s,x}(\Gamma^{c})\leq\mathbb{P}^{s,x}\left(\max_{\substack{t\in\mathcal{I}_{i} \\ i=1,\ldots,M_{1}}}\left|X_{t}^{(\rho,n)}-X_{s_{i}}^{(\rho,n)}\right|\geq\frac{\kappa_{1}}{2}\right)\leq M_{1}^{-1}\kappa_{1}^{-4}D_{1}T^{2}.
\end{align}
Hence, by separating the expectation in $\Gamma$ and $\Gamma^{c}$, the second term in \eqref{eq:DynkinWRhon2} can be estimated as
\begin{align*} \mathbb{E}^{s,x}\!\!\left(\!\int_{s}^{\tau\wedge\tau_{\rho,n}^{(\delta)}}\!\!\left(\mathscr{A}_{\rho,\widehat{U}}^{(n)}\mathcal{W}^{(n)}_{\rho}\!\!\left(t,X_{t}^{(\rho,n)}\right)\!+\!\mathcal{L}_{\rho}\!\!\left(X_{(\rho,n)}^{t},\widehat{U}_{t}\right)\!\right)\!dt\!\right)\!\leq\!\frac{\varepsilon}{2}\!+\!\frac{D_{1}T^{2}}{M_{1}\kappa_{1}^{4}}\!\left\|\mathscr{A}_{\rho,u}^{(n)}\mathcal{W}^{(n)}_{\rho}\!+\!\mathcal{L}_{\rho}\right\|_{L^{\infty}(\overline{Q}_{\rho-\delta,T}\times\mathscr{U})}.
\end{align*}
Therefore, for fixed $\rho$, $\delta>0$, when $M_{1}>0$ is large enough, for any reference stochastic system $\mu\in\widetilde{\mathscr{P}}_{[s,T]}$, there exists $\widehat{U}\in\mathcal{V}_{\mu}[s,T]$, such that for any $(\mathscr{F}_{t})_{t\in[s,T]}$-stopping time $\tau$,
\begin{align}\label{eq:DPRhonDelta1}
\mathcal{W}^{(n)}_{\rho}(s,x)+\varepsilon\geq\mathbb{E}^{s,x}\left(\int_{s}^{\tau\wedge\tau_{\rho,n}^{(\delta)}}\mathcal{L}_{\rho}\!\left(X_{t}^{(\rho,n)},\widehat{U}_{t}\right)dt+\mathcal{W}^{(n)}_{\rho}\!\left(\tau\wedge\tau_{\rho,n}^{(\delta)},X_{\tau\wedge\tau_{\rho,n}^{(\delta)}}^{(\rho,n)}\right)\right).
\end{align}
Also, by \eqref{eq:RhonHJB2} and \eqref{eq:DynkinWRhon1}, for any $\mu\in\widetilde{\mathscr{P}}_{[s,T]}$, $U\in\mathcal{V}_{\mu}[s,T]$, and any $(\mathscr{F}_{t})_{t\in[s,T]}$-stopping time $\tau$,
\begin{align}\label{eq:DPRhonDelta2}
\mathcal{W}^{(n)}_{\rho}(s,x)\leq\mathbb{E}_{s,x}\left(\int_{s}^{\tau\wedge\tau_{\rho,n}^{(\delta)}}\mathcal{L}_{\rho}\!\left(X_{t}^{(\rho,n)},U_{t}\right)dt+\mathcal{W}^{(n)}_{\rho}\!\left(\tau\wedge\tau_{\rho,n}^{(\delta)},X_{\tau\wedge\tau_{\rho,n}^{(\delta)}}^{(\rho,n)}\right)\right).
\end{align}

We will take $\delta\rightarrow 0$ in both \eqref{eq:DPRhonDelta1} and \eqref{eq:DPRhonDelta2}. Since $\mathcal{W}^{(n)}_{\rho}$ is uniformly continuous on $\overline{Q}_{\rho,T}$, there exists $\kappa_{2}>0$ so that, for any $(t,y),(t',y')\in\overline{Q}_{\rho,T}$ with $|t-t'|\leq\kappa_{2}$ and $|y-y'|\leq\kappa_{2}$,
\begin{align*}
\left|\mathcal{W}^{(n)}_{\rho}(t,y)-\mathcal{W}^{(n)}_{\rho}(t',y')\right|\leq\varepsilon.
\end{align*}
Also, for any $\delta<\tilde{\kappa}$, the event $\{\tau_{\rho,n}-\tau_{\rho,n}^{(\delta)}>\tilde{\kappa}\}$ occurs only when $X^{(\rho,n)}$ hits $\partial\mathscr{O}_{\rho-\delta}$ before $T-\delta$ (otherwise, $\tau_{\rho,n}^{(\delta)}=T-\delta$ and $\tau_{\rho,n}\leq T$). By conditioning on $\tau_{\rho,n}^{(\delta)}$, we have
\begin{align*}
\mathbb{P}^{s,x}\left(\tau_{\rho,n}-\tau_{\rho,n}^{(\delta)}>\tilde{\kappa}\right)\leq\tilde{\kappa}^{-1}\psi(\delta).
\end{align*}
Together with Assumption \ref{assump:SmallBallCondRhon}, for $\delta<\kappa_{2}$, we have
\begin{align}
&\mathbb{E}^{s,x}\left(\left|\mathcal{W}_{\rho}^{(n)}\left(\tau\wedge\tau_{\rho,n}^{(\delta)},X_{\tau\wedge\tau_{\rho,n}^{(\delta)}}^{(\rho,n)}\right)-\mathcal{W}_{\rho}^{(n)}\left(\tau\wedge\tau_{\rho,n},X_{\tau\wedge\tau_{\rho,n}}^{(\rho,n)}\right)\right|\right)\nonumber\\
&\quad\leq\varepsilon+2\max_{(t,y)\in\overline{Q}_{\rho,T}}\left|\mathcal{W}_{\rho}^{(n)}(t,y)\right|\left(\mathbb{P}^{s,x}\left(\tau_{\rho,n}-\tau_{\rho,n}^{(\delta)}>\kappa_{2}\right)+\mathbb{P}^{s,x}\left(\left|X_{\tau\wedge\tau_{\rho,n}}^{(\rho,n)}-X_{\tau\wedge\tau_{\rho,n}^{(\delta)}}^{(\rho,n)}\right|>\kappa_{2}\right)\right)\nonumber\\
&\quad\leq\varepsilon+2\max_{(t,y)\in\overline{Q}_{\rho,T}}\left|\mathcal{W}_{\rho}^{(n)}(t,y)\right|\left(\mathbb{P}^{s,x}\left(\tau_{\rho,n}-\tau_{\rho,n}^{(\delta)}>\kappa_{2}\right)+\mathbb{P}^{s,x}\left(\tau_{\rho,n}-\tau_{\rho,n}^{(\delta)}>2\delta\right)\right.\nonumber\\
&\qquad\qquad\qquad\qquad\qquad\qquad\qquad\left.+\mathbb{P}^{s,x}\left(\max_{t,t'\in[s,T],\,|t-t'|<2\delta}\left|X_{t}^{(\rho,n)}-X_{t'}^{(\rho,n)}\right|>\kappa_{2}\right)\right)\nonumber\\
\label{eq:EstDiffWRhonDelta} &\quad\leq\varepsilon+2\max_{(t,y)\in\overline{Q}_{\rho,T}}\left|\mathcal{W}_{\rho}^{(n)}(t,y)\right|\left(\kappa_{2}^{-1}\psi(\delta)+(2\delta)^{-1}\psi(2\delta)+2\kappa_{2}^{-1}D_{1}\delta\right),
\end{align}
where we used again (D.12) in~\cite{FlemingSoner:2006} to estimate the last probability above, and where $D_{1}>0$ is as in \eqref{eq:EstGammaComp}, and only depends on $\rho$, $T$ and $K$. Moreover, for any $U\in\mathcal{V}_{\mu}[s,T]$,
\begin{align}\label{eq:EstDiffLRhonDelta}
\mathbb{E}^{s,x}\!\left(\int_{\tau\wedge\tau_{\rho,n}^{(\delta)}}^{\tau\wedge\tau_{\rho,n}}\!\left|\mathcal{L}_{\rho}\!\left(X_{t}^{(\rho,n)},U_{t}\right)\right|dt\right)\leq\kappa\mathscr{L}^{\gamma}\rho\,\mathbb{E}^{s,x}\!\left(\tau_{\rho,n}-\tau_{\rho,n}^{(\delta)}\right)\leq\kappa\mathscr{L}^{\gamma}\rho\left(T\delta^{-1}\psi(\delta)+\delta\right).
\end{align}
Noting that $\delta^{-1}\psi(\delta)\rightarrow\psi'_{+}(0)=0$, as $\delta\rightarrow 0$, by combining \eqref{eq:EstDiffWRhonDelta} and \eqref{eq:EstDiffLRhonDelta}, for any $(s,x)\in Q_{\rho,T}$, $\mu\in\widetilde{\mathscr{P}}_{[s,T]}$, $U\in\mathcal{V}_{\mu}[s,T]$, and any $(\mathscr{F}_{t})_{t\in[s,T]}$-stopping time $\tau$, we have shown that
\begin{align}\label{eq:LimDiffWLRhonDelta}
\lim_{\delta\rightarrow
0}\sup_{\tau}\mathbb{E}^{s,x}\!\!\left(\!\left|\mathcal{W}^{(n)}_{\rho}\!\!\left(\!\tau\!\wedge\!\tau_{\rho,n}^{(\delta)},\!X_{\tau\wedge\tau_{\rho,n}^{(\delta)}}^{(\rho,n)}\!\right)\!-\!\mathcal{W}^{(n)}_{\rho}\!\!\left(\!\tau\!\wedge\!\tau_{\rho,n},\!X_{\tau\wedge\tau_{\rho,n}}^{(\rho,n)}\!\right)\!\right|\!+\!\!\int_{\tau\wedge\tau_{\rho,n}^{(\delta)}}^{\tau\wedge\tau_{\rho,n}}\!\!\left|\mathcal{L}_{\rho}\!\!\left(\!X_{t}^{(\rho,n)}\!,U_{t}\!\right)\right|\!dt\!\right)\!\!=\!0.
\end{align}
Hence, letting $\delta\rightarrow 0$ in \eqref{eq:DPRhonDelta2}, for any $(s,x)\in Q_{\rho,T}$, $\mu\in\widetilde{\mathscr{P}}_{[s,T]}$, $U\in\mathcal{V}_{\mu}[s,T]$, and any $(\mathscr{F}_{t})_{t\in[s,T]}$-stopping time $\tau$,
\begin{align}\label{eq:DPRhon1}
\mathcal{W}^{(n)}_{\rho}(s,x)\leq\mathbb{E}_{s,x}\left(\!\int_{s}^{\tau\wedge\tau_{\rho,n}}\mathcal{L}_{\rho}\!\left(X_{t}^{(\rho,n)},U_{t}\right)dt+\mathcal{W}^{(n)}_{\rho}\left(\tau\wedge\tau_{\rho,n},X_{\tau\wedge\tau_{\rho,n}}^{(\rho,n)}\right)\right).
\end{align}
Moreover, by \eqref{eq:LimDiffWLRhonDelta}, for any $\varepsilon>0$, $\mu\in\widetilde{\mathscr{P}}_{[s,T]}$, $U\in\mathcal{V}_{\mu}[s,T]$, there exists $\delta_{1}>0$, such that for any $\delta\in(0,\delta_{1}]$ and any $(\mathscr{F}_{t})_{t\in[s,T]}$-stopping time $\tau$,
\begin{align*}
\mathbb{E}^{s,x}\left(\left|\mathcal{W}^{(n)}_{\rho}\!\left(\tau\wedge\tau_{\rho,n}^{(\delta)},X_{\tau\wedge\tau_{\rho,n}^{(\delta)}}^{(\rho,n)}\right)-\mathcal{W}^{(n)}_{\rho}\!\left(\tau\wedge\tau_{\rho,n},X_{\tau\wedge\tau_{\rho,n}}^{(\rho,n)}\right)\right|+\int_{\tau\wedge\tau_{\rho,n}^{(\delta)}}^{\tau\wedge\tau_{\rho,n}}\left|\mathcal{L}_{\rho}\!\left(X_{t}^{(\rho,n)},U_{t}\right)\right|dt\right)\leq\frac{\varepsilon}{2}.
\end{align*}
Combining the above together with \eqref{eq:DPRhonDelta1} (with $\varepsilon$ replaced by $\varepsilon/2$, and choosing $\delta_{2}\in(0,\delta_{1})$ so that $(s,x)\in\overline{Q}_{\rho-\delta_{2},T}$), we obtain that, for any $\mu\in\widetilde{\mathscr{P}}_{[s,T]}$, there exists $\widehat{U}\in\mathcal{V}_{\mu}[s,T]$, such that for any $(\mathscr{F}_{t})_{t\in[s,T]}$-stopping time $\tau$,
\begin{align}\label{eq:DPRhon2} \mathcal{W}^{(n)}_{\rho}(s,x)+\varepsilon\geq\mathbb{E}^{s,x}\left(\int_{s}^{\tau\wedge\tau_{\rho,n}}\mathcal{L}_{\rho}\!\left(X_{t}^{(\rho,n)},\widehat{U}_{t}\right)dt+\mathcal{W}^{(n)}_{\rho}\!\left(\tau\wedge\tau_{\rho,n},X_{\tau\wedge\tau_{\rho,n}}^{(\rho,n)}\right)\right).
\end{align}
Note that both \eqref{eq:DPRhon1} and \eqref{eq:DPRhon2} are trivially true for $(s,x)\in\partial^{\ast}Q_{\rho}$, and hence for all $(s,x)\in\overline{Q}_{\rho,T}$. In particular, letting $\theta\equiv T$, then for any $\mu\in\widetilde{\mathscr{P}}_{[s,T]}$,
\begin{align}\label{eq:DynProRhon}
\mathcal{W}^{(n)}_{\rho}(s,x)=\inf_{U\in\mathcal{V}_{\mu}[s,T]}\mathbb{E}^{s,x}\left(\int_{s}^{\tau_{\rho,n}}\mathcal{L}_{\rho}\left(X_{t}^{(\rho,n)},U_{t}\right)dt-P_{\tau_{\rho,n}}^{(\rho,n)}Z_{\tau_{\rho,n}}^{(\rho,n)}\right)=V_{\rho,\mu}^{(n)}(s,y),
\end{align}
which immediately implies that for all $(s,x)\in\overline{Q}_{\rho,T}$ and $\mu\in\widetilde{\mathscr{P}}_{[s,T]}$,
\begin{align*}
\mathcal{W}^{(n)}_{\rho}(s,x)=V_{\rho,\mu}^{(n)}(s,x)=V_{\rho}^{(n)}(s,x)\in C(\overline{Q}_{\rho,T}).
\end{align*}

\noindent
\textbf{Step 2.} Fix any $\rho>0$. For any $(s,x)\in\overline{Q}_{\rho,T}$, $\nu\in\mathscr{P}_{[s,T]}$, and $U\in\mathcal{U}_{\nu}[s,T]$, consider the SDE
\begin{align*}
dX_{t}^{(\rho)}=\vec{f_{\rho}}\left(t,X_{t}^{(\rho)},U_{t}\right)dt+\vec{\sigma}_{\rho}\left(t,X_{t}^{(\rho)},U_{t}\right)dW_{t},
\end{align*}
with initial condition $X_{s}^{(\rho)}=x$, $\mathbb{P}$-a.$\,$s., and the associated stochastic control problem
\begin{align}\label{eq:RhoCostFunt}
J_{\rho,\nu}(s,x;U)&:=\mathbb{E}^{s,x}\left(\int_{s}^{\tau_{\rho}}\mathcal{L}_{\rho}\left(X_{t}^{(\rho)},U_{t}\right)dt-P_{\tau_{\rho}}^{(\rho)}Z_{\tau_{\rho}}^{(\rho)}\right),\\
\label{eq:RhoValFunt1} V_{\rho,\nu}(s,x)&:=\inf_{U\in\mathcal{U}_{\nu}[s,T]}J_{\rho,\nu}(s,x;U),\\
\label{eq:RhoValFunt2} V_{\rho}(s,x)&:=\inf_{\nu\in\mathscr{P}_{[s,T]}}V_{\rho,\nu}(s,x),
\end{align}
where we set $X^{(\rho)}=(P^{(\rho)},Z^{(\rho)},\Theta^{(\rho)})$, and where
\begin{align*}
\tau_{\rho}=\tau_{\rho}(s,x):=\inf\left\{t\geq s:\,\,X_{t}^{(\rho)}\notin\mathscr{O}_{\rho},\,\,X_{s}^{(\rho)}=x\right\}\wedge T.
\end{align*}
By Remark \ref{rem:EquiStoSym}, we can build a one-to-one correspondence between all $\mu\in\widetilde{\mathscr{P}}_{[s,T]}$ and $\nu\in\mathscr{P}_{[s,T]}$. Hence, we can define the expectation in $J_{\rho,\nu}(s,x;U)$ on the same six-tuple stochastic system $\mu\in\widetilde{\mathscr{P}}_{[s,T]}$ as $J_{\rho,\nu}^{(n)}(s,x;U)$, and choose the control policy $U\in\mathcal{V}_{\mu}[s,T]$ for both $J_{\rho,\nu}(s,x;U)$ and $J_{\rho,\nu}(s,x;U)$. We will prove the uniform convergence of $J_{\rho,\mu}^{(n)}(s,x;U)$ towards $J_{\rho,\mu}(s,x;U)$, as $n\rightarrow\infty$, with respect to all $(s,x)\in\overline{Q}_{\rho,T}$, $\mu\in\widetilde{\mathscr{P}}_{[s,T]}$ and $U\in\mathcal{V}_{\mu}[s,T]$.

To see this, for any $(s,x)\in\overline{Q}_{\rho,T}$, $\mu\in\widetilde{\mathscr{P}}_{[s,T]}$ and $U\in\mathcal{V}_{\mu}[s,T]$,
\begin{align}
&\left|J_{\rho,\mu}(s,x;U)-J_{\rho,\mu}^{(n)}(s,x;U)\right|\nonumber\\
&\quad\leq\mathbb{E}^{s,x}\left(\left|\int_{s}^{\tau_{\rho}}\mathcal{L}_{\rho}\!\left(X_{t}^{(\rho)},U_{t}\right)dt-\int_{s}^{\tau_{\rho,n}}\mathcal{L}_{\rho}\!\left(X_{t}^{(\rho,n)},U_{t}\right)dt\right|\right)+\mathbb{E}^{s,x}\left(\left|P^{(\rho)}_{\tau_{\rho}}Z^{(\rho)}_{\tau_{\rho}}-P^{(\rho,n)}_{\tau_{\rho,n}}Z^{(\rho,n)}_{\tau_{\rho,n}}\right|\right)\nonumber\\
&\quad\leq\mathbb{E}^{s,x}\!\!\left(\!\int_{s}^{\tau_{\rho}\wedge\tau_{\rho,n}}\!\!\left|\mathcal{L}_{\rho}\!\left(\!X_{t}^{(\rho)}\!,U_{t}\!\right)\!-\!\mathcal{L}_{\rho}\!\left(\!X_{t}^{(\rho,n)}\!,U_{t}\!\right)\right|\!dt\!\right)\!\!+\!\kappa\mathscr{L}^{\gamma}\!\rho\mathbb{E}^{s,x}\!\!\left(\left|\tau_{\rho,n}\!-\!\tau_{\rho}\right|\right)\!+\!2\rho\mathbb{E}^{s,x}\!\!\left(\left|\!X_{\tau_{\rho}}^{(\rho)}\!-\!X_{\tau_{\rho,n}}^{(\rho,n)}\!\right|\right)\nonumber\\
\label{eq:DecompJRhoJRhon} &\quad\leq\kappa\mathscr{L}^{\gamma}\!\!\int_{s}^{T}\!\mathbb{E}^{s,x}\!\left(\left|X_{t}^{(\rho)}\!-\!X_{t}^{(\rho,n)}\right|\right)dt+\kappa\mathscr{L}^{\gamma}\rho\,\mathbb{E}^{s,x}\!\left(\left|\tau_{\rho,n}\!-\!\tau_{\rho}\right|\right)+2\rho\mathbb{E}^{s,x}\!\left(\left|X_{\tau_{\rho}}^{(\rho)}\!-\!X_{\tau_{\rho,n}}^{(\rho,n)}\right|\right).
\end{align}
For the first term in \eqref{eq:DecompJRhoJRhon}, by (D.9) in~\cite{FlemingSoner:2006}, there exists a constant $D_{2}>0$, depending only on $\rho$, $T$ and $K$, such that for any $t\in[s,T]$,
\begin{align}\label{eq:EstDiffJRhoJRhon1}
\mathbb{E}^{s,x}\left(\left|X_{t}^{(\rho)}-X_{t}^{(\rho,n)}\right|\right)\leq\mathbb{E}^{s,x}\left(\max_{t\in[s,T]}\left|X_{t}^{(\rho)}-X_{t}^{(\rho,n)}\right|\right)\leq D_{2}\,\epsilon^{n},
\end{align}
since $X^{(\rho)}$ and $X^{(\rho,n)}$ only differ in a diffusion term of $\epsilon^{n}$. To estimate the second term in \eqref{eq:DecompJRhoJRhon}, for any $\varepsilon>0$, first pick $\delta=\delta(\varepsilon)>0$ small enough so that $\psi(\delta)\leq\varepsilon^{2}$. By Assumption \ref{assump:SmallBallCond}, Assumption \ref{assump:SmallBallCondRhon} and \eqref{eq:EstDiffJRhoJRhon1},
\begin{align}
\mathbb{E}^{s,x}\left(\left|\tau_{\rho,n}-\tau_{\rho}\right|\right)&\leq\varepsilon+2T\,\mathbb{P}^{s,x}\left(\left|\tau_{\rho,n}-\tau_{\rho}\right|>\varepsilon,\,\max_{t\in[s,T]}\left|X_{t}^{(\rho)}-X_{t}^{(\rho,n)}\right|\leq\delta\right)+2TD_{2}\delta^{-1}\epsilon^{n}\nonumber\\
\label{eq:EstDiffJRhoJRhon2} &\leq\varepsilon+2T\varepsilon^{-1}\psi(\delta)+2TD_{2}\delta^{-1}\epsilon^{n}\leq(1+2T)\varepsilon+2TD_{2}\delta^{-1}\epsilon^{n}.
\end{align}
Moreover, for the last term in \eqref{eq:DecompJRhoJRhon}, by \eqref{eq:EstDiffJRhoJRhon1}, and for any $\varepsilon>0$,
\begin{align*}
\mathbb{E}^{s,x}\left(\left|X_{\tau_{\rho}}^{(\rho)}-X_{\tau_{\rho,n}}^{(\rho,n)}\right|\right)&\leq\mathbb{E}^{s,x}\left(\left|X_{\tau_{\rho}}^{(\rho)}-X_{\tau_{\rho,n}}^{(\rho)}\right|\right)+\mathbb{E}^{s,x}\left(\left|X_{\tau_{\rho,n}}^{(\rho)}-X_{\tau_{\rho,n}}^{(\rho,n)}\right|\right)\\
&\leq\varepsilon+2\rho\,\mathbb{P}^{s,x}\left(\left|X_{\tau_{\rho}}^{(\rho)}-X_{\tau_{\rho,n}}^{(\rho)}\right|>\varepsilon\right)+D_{2}\epsilon^{n}\\
&\leq\varepsilon\!+\!2\rho\!\left(\mathbb{P}^{s,x}\!\left(\left|\tau_{\rho,n}\!-\!\tau_{\rho}\right|\!>\!\varepsilon^{3}\right)\!+\!\mathbb{P}^{s,x}\!\left(\left|\!X_{\tau_{\rho}}^{(\rho)}\!-\!X_{\tau_{\rho,n}}^{(\rho)}\!\right|\!>\!\varepsilon,\left|\tau_{\rho,n}\!-\!\tau_{\rho}\right|\!\leq\!\varepsilon^{3}\right)\!\right)\!+\!D_{2}\epsilon^{n}.
\end{align*}
In a similar fashion to \eqref{eq:EstDiffJRhoJRhon2}, but with $\delta=\delta(\varepsilon)>0$ so that $\psi(\delta)\leq\varepsilon^{4}$, we have
\begin{align*}
\mathbb{P}^{s,x}\left(\left|\tau_{\rho,n}-\tau_{\rho}\right|>\varepsilon^{3}\right)\leq\varepsilon+D_{2}\delta^{-1}\epsilon^{n}.
\end{align*}
Moreover, by conditioning on $\mathscr{F}_{\tau_{\rho}\wedge\tau_{\rho,n}}$, and using the strong Markov property of $X^{(\rho)}$ as well as (D.12) in~\cite{FlemingSoner:2006},
\begin{align*}
&\mathbb{P}^{s,x}\left(\left|X_{\tau_{\rho}}^{(\rho)}-X_{\tau_{\rho,n}}^{(\rho)}\right|>\varepsilon,\,\left|\tau_{\rho,n}-\tau_{\rho}\right|\leq\varepsilon^{3}\right)\\
&\quad=\int_{\overline{Q}_{\rho,T}}\mathbb{P}^{s,x}\!\left(\left.\left|X_{\tau_{\rho}}^{(\rho)}\!-\!X_{\tau_{\rho,n}}^{(\rho)}\right|\!>\!\varepsilon,\left|\tau_{\rho,n}\!-\!\tau_{\rho}\right|\!\leq\!\varepsilon^{3}\,\right|\left(\tau_{\rho}\!\wedge\!\tau_{\rho,n},X_{\tau_{\rho}\wedge\tau_{\rho,n}}^{(\rho)}\right)\!=\!(t,y)\right)dF_{\tau_{\rho}\wedge\tau_{\rho,n},X_{\tau_{\rho}\wedge\tau_{\rho,n}}^{(\rho)}}\!\!(t,y)\\
&\quad\leq\int_{\overline{Q}_{\rho,T}}\mathbb{P}^{t,y}\left(\max_{u\in[t,t+\varepsilon^{3}]}\left|X_{u}^{(\rho)}-y\right|>\varepsilon\right)dF_{\tau_{\rho}\wedge\tau_{\rho,n},X_{\tau_{\rho}\wedge\tau_{\rho,n}}^{(\rho)}}\!\!(t,y)\leq D_{3}\varepsilon,
\end{align*}
where $D_{3}>0$ is a constant depending on $\rho$, $T$ and $K$. Hence,
\begin{align}\label{eq:EstDiffJRhoJRhon3}
\mathbb{E}^{s,x}\left(\left|X_{\tau_{\rho}}^{(\rho)}-X_{\tau_{\rho,n}}^{(\rho,n)}\right|\right)\leq(1+2\rho+D_{3})\varepsilon+\left(2\rho\delta^{-1}+1\right)D_{2}\epsilon^{n}.
\end{align}

Combining \eqref{eq:DecompJRhoJRhon}-\eqref{eq:EstDiffJRhoJRhon3} shows that $J_{\rho,\mu}^{(n)}(s,x;U)$ converges, as $n\rightarrow\infty$, to $J_{\rho,\mu}(s,x;U)$ uniformly for all $(s,x)\in\overline{Q}_{\rho,T}$, $\mu\in\widetilde{\mathscr{P}}_{[s,T]}$ and $U\in\mathcal{V}_{\mu}[s,T]$. It then immediately follows that as $n\rightarrow\infty$, $V_{\rho,\mu}^{(n)}(s,x)\rightarrow V_{\rho,\mu}(s,x)$ uniformly in $(s,x)\in\overline{Q}_{\rho,T}$ and $\mu\in\widetilde{\mathscr{P}}_{[s,T]}$, and that $V_{\rho}^{(n)}(s,x)\rightarrow V_{\rho}(s,x)$ uniformly in $(s,x)\in\overline{Q}_{\rho,T}$. Therefore, $V_{\rho}\in C(\overline{Q}_{\rho,T})$ and $V_{\rho}(s,x)=V_{\rho,\mu}(s,x)$ for all $(s,x)\in\overline{Q}_{\rho,T}$ and $\mu\in\widetilde{\mathscr{P}}_{[s,T]}$. By the one-to-one correspondence between the collection of six-tuple stochastic systems $\widetilde{\mathscr{P}}_{[s,T]}$ and the collection of five-tuple stochastic systems $\mathscr{P}_{[s,T]}$, we conclude that
$V_{\rho}(s,x)=V_{\rho,\nu}(s,x)$ for all $(s,x)\in\overline{Q}_{\rho,T}$ and $\nu\in\mathscr{P}_{[s,T]}$.

\medskip
\noindent
\textbf{Step 3.} We now consider the stochastic control problem \eqref{eq:WeakCostFunt}-\eqref{eq:WeakValFunt2}. In similarity to Step 2, we will prove the uniform convergence of $J_{\rho,\nu}(s,x;U)$ towards $J_{\nu}(s,x;U)$, as $\rho\rightarrow\infty$, for all $(s,x)\in[0,T]\times\mathcal{C}$, $\nu\in\mathscr{P}_{[s,T]}$ and $U\in\mathcal{V}_{\nu}[s,T]$, where $\mathcal{C}$ is an arbitrary compact subset of $\overline{\mathscr{O}}\setminus\{z=0\}$.

For such a $\mathcal{C}\subseteq\overline{\mathscr{O}}\setminus\{z=0\}$, we can find $\rho>0$ large enough so that $\mathcal{C}\subseteq\mathscr{O}_{\rho}$. For any $(s,x)\in[0,T]\times\mathcal{C}$, $\nu\in\mathscr{P}_{[s,T]}$ and $U\in\mathcal{V}_{\nu}[s,T]$, noticing that $\tau_{\rho}\leq\tau_{\mathscr{O}}$, and that the trajectories of $X^{(\rho)}$ and $X$ are identical up to time $\tau_{\rho}$, we have
\begin{align}
&\left|J_{\rho,\nu}(s,x;U)-J_{\nu}(s,x;U)\right|\leq\mathbb{E}^{s,x}\left(\int_{\tau_{\rho}}^{\tau_{\mathscr{O}}}\kappa L_{t}^{\gamma}Z_{t}\,dt\right)+\mathbb{E}^{s,x}\left(\left|P_{\tau_{\mathscr{O}}}Z_{\tau_{\mathscr{O}}}-P_{\tau_{\rho}}Z_{\tau_{\rho}}\right|\right)\nonumber\\
\label{eq:DecompJRhoJ} &\quad\leq\kappa\mathscr{L}^{\gamma}\mathbb{E}^{s,x}\!\left(\left|\tau_{\mathscr{O}}\!-\!\tau_{\rho}\right|\left(\sup_{t\in[s,T]}Z_{t}\right)\right)+\mathbb{E}^{s,x}\!\left(\left|P_{\tau_{\mathscr{O}}}\right|\left|Z_{\tau_{\mathscr{O}}}\!-\!Z_{\tau_{\rho}}\right|\right)+\mathbb{E}^{s,x}\!\left(\left|Z_{\tau_{\rho}}\right|\left|P_{\tau_{\mathscr{O}}}\!-\!P_{\tau_{\rho}}\right|\right).
\end{align}
For the first expectation in \eqref{eq:DecompJRhoJ}, we first have
\begin{align*}
\mathbb{E}^{s,x}\!\left(\left|\tau_{\mathscr{O}}-\tau_{\rho}\right|\left(\sup_{t\in[s,T]}Z_{t}\right)\right)\leq\sqrt{\mathbb{E}^{s,x}\!\left(\left|\tau_{\mathscr{O}}-\tau_{\rho}\right|^{2}{\bf 1}_{\{\tau_{\rho}<\tau_{\mathscr{O}}\}}\right)\mathbb{E}^{s,x}\!\left(\sup_{t\in[s,T]}Z_{t}^{2}\right)}.
\end{align*}
By the construction of $\mathscr{O}_{\rho}$, for any $\varepsilon>0$,
\begin{align*}
\mathbb{E}^{s,x}\!\left(\left|\tau_{\mathscr{O}}-\tau_{\rho}\right|^{2}{\bf 1}_{\{\tau_{\rho}<\tau_{\mathscr{O}}\}}\right)&\leq 2T\,\mathbb{P}^{s,x}\!\left(\sup_{t\in[s,T]}\left|P_{t}\right|\geq\rho\right)+2T\,\mathbb{P}^{s,x}\!\left(\sup_{t\in[s,T]}Z_{t}\geq\rho\right)\\ &\quad\,+2T\,\mathbb{P}^{s,x}\!\left(\text{dist}\!\left(X_{\tau_{\rho}},\partial\mathscr{O}\right)\leq\frac{1}{\rho},\,\tau_{\mathscr{O}}-\tau_{\rho}>\varepsilon\right)+\varepsilon.
\end{align*}
By Assumption \ref{assump:Coefficients}$-$(i) and (D.7) in~\cite{FlemingSoner:2006}, there exists $D_{4}>0$, depending on $K$ and $T$, so that
\begin{align}\label{eq:EstMaxTailP}
\mathbb{P}^{s,x}\!\left(\sup_{t\in[s,T]}\left|P_{t}\right|\geq\rho\right)\leq\frac{1}{\rho}\,\mathbb{E}^{s,x}\!\left(\sup_{t\in[s,T]}\left|P_{t}\right|\right)\leq\frac{D_{4}}{\rho}\left(1+\max_{x\in\mathcal{C}}|x|\right),
\end{align}
and by Doob's martingale inequality (cf.~\cite[Theorem 1.3.8-(i)]{KaratzasShreve:1991}),
\begin{align}\label{eq:EstMaxTailZ}
\mathbb{P}^{s,x}\!\left(\sup_{t\in[s,T]}Z_{t}\geq\rho\right)\leq\frac{1}{\rho}\,\mathbb{E}(Z_{T})=\frac{z}{\rho}\leq\frac{1}{\rho}\max_{x\in\mathcal{C}}|x|.
\end{align}
Moreover, by conditioning on $\mathscr{F}_{\tau_{\rho}}$ and using Assumption \ref{assump:SmallBallCond},
\begin{align}\label{eq:EstProbDiffTauTauRho}
\mathbb{P}^{s,x}\!\left(\text{dist}\!\left(X_{\tau_{\rho}},\partial\mathscr{O}\right)\leq\frac{1}{\rho},\,\tau_{\mathscr{O}}-\tau_{\rho}>\varepsilon\right)\leq\varepsilon^{-1}\psi(\rho^{-1}).
\end{align}
Finally, by Doob's martingale inequality (cf.~\cite[Theorem 1.3.8$-$(iv)]{KaratzasShreve:1991}),
\begin{align}
\mathbb{E}^{s,x}\!\left(\sup_{t\in[s,T]}Z_{t}^{2}\right)\leq 4\mathbb{E}^{s,x}\!\left(Z_{T}^{2}\right)&=4\mathbb{E}^{s,x}\!\left(z^{2}e^{2\int_{s}^{T}\left(\Theta_{u}+A\,\Xi_{u}^{\alpha}L_{u}^{\beta}-C_{u}\right)\varrho^{-1}dW_{u}-\int_{s}^{T}\left(\Theta_{u}+A\,\Xi_{u}^{\alpha}L_{u}^{\beta}-C_{u}\right)^{2}\varrho^{-2}du}\right)\nonumber\\
\label{eq:EstSupZtSqu} &\leq 4\max_{x\in\mathcal{C}}|x|^{2}\exp\left(\left(H+A\mathscr{N}^{\alpha}\mathscr{L}^{\beta}+\mathscr{C}\right)^{2}\varrho^{-2}T\right).
\end{align}
Hence, the first expectation in \eqref{eq:DecompJRhoJ} can be estimated by
\begin{align}\label{eq:EstDiffJRhoJ1}
\mathbb{E}^{s,x}\!\left(\left|\tau_{\mathscr{O}}-\tau_{\rho}\right|\left(\sup_{t\in[s,T]}Z_{t}\right)\right)\leq\sqrt{24D_{4}T\max_{x\in\mathcal{C}}|x|^{3}e^{\left(H+A\mathscr{N}^{\alpha}\mathscr{L}^{\beta}+\mathscr{C}\right)^{2}\varrho^{-2}T}\!\left(\frac{1}{\rho}\!+\!\frac{1}{\varepsilon}\psi(\rho^{-1})\!+\!\varepsilon\!\right)}.
\end{align}

The second and the third expectations in \eqref{eq:DecompJRhoJ} can be analyzed in the same way, and so only the estimation for the second expectation in \eqref{eq:DecompJRhoJ} is presented. For any $\varepsilon>0$,
\begin{align*}
\mathbb{E}^{s,x}\!\left(\left|P_{\tau_{\mathscr{O}}}\right|\left|Z_{\tau_{\mathscr{O}}}\!-Z_{\tau_{\rho}}\right|\right)=\mathbb{E}^{s,x}\!\left(\left|P_{\tau_{\mathscr{O}}}\right|\left|Z_{\tau_{\mathscr{O}}}\!-Z_{\tau_{\rho}}\right|{\bf 1}_{\{\tau_{\mathscr{O}}-\tau_{\rho}\leq\varepsilon\}}\right)+\mathbb{E}^{s,x}\!\left(\left|P_{\tau_{\mathscr{O}}}\right|\left|Z_{\tau_{\mathscr{O}}}\!-Z_{\tau_{\rho}}\right|{\bf 1}_{\{\tau_{\mathscr{O}}-\tau_{\rho}>\varepsilon\}}\right).
\end{align*}
By \eqref{eq:EstSupZtSqu} and (D.7) in~\cite{FlemingSoner:2006},
\begin{align*}
\mathbb{E}^{s,x}\!\left(\left|P_{\tau_{\mathscr{O}}}\right|\!\left|Z_{\tau_{\mathscr{O}}}\!\!-\!Z_{\tau_{\rho}}\right|\!{\bf 1}_{\{\tau_{\mathscr{O}}-\tau_{\rho}>\varepsilon\}}\right)&\leq\left(\!\mathbb{E}^{s,x}\!\!\left(\sup_{t\in[s,T]}\!\left|P_{t}\right|^{4}\!\right)\!\right)^{1/4}\!\!\!\!\!\left(\mathbb{P}^{s,x}\!\left(\tau_{\mathscr{O}}\!-\!\tau_{\rho}>\varepsilon\right)\right)^{1/4}\!\left(\!\mathbb{E}^{s,x}\!\left(\sup_{t\in[s,T]}\!Z_{t}^{2}\!\right)\!\right)^{1/2}\\
&\leq D_{5}\max_{x\in\mathcal{C}}|x|^{2}\exp\!\left(\!\frac{T}{2\varrho^{2}}\!\left(\!H\!+\!A\mathscr{N}^{\alpha}\!\mathscr{L}^{\beta}\!+\!\mathscr{C}\!\right)^{2}\right)\!\left(\mathbb{P}^{s,x}\!\left(\tau_{\mathscr{O}}\!-\!\tau_{\rho}\!>\!\varepsilon\right)\right)^{1/4},
\end{align*}
where $D_{5}>0$ is a constant depending on $K$ and $T$. Moreover, by \eqref{eq:EstMaxTailP}-\eqref{eq:EstProbDiffTauTauRho},
\begin{align}
\mathbb{P}^{s,x}\!\left(\tau_{\mathscr{O}}\!-\!\tau_{\rho}\!>\!\varepsilon\right)&\leq\mathbb{P}^{s,x}\!\left(\sup_{t\in[s,T]}\!\left|P_{t}\right|\!\geq\!\rho\!\right)+\mathbb{P}^{s,x}\!\left(\sup_{t\in[s,T]}\!Z_{t}\!\geq\!\rho\!\right)+\mathbb{P}^{s,x}\!\left(\text{dist}\!\left(X_{\tau_{\rho}},\partial\mathscr{O}\right)\!\leq\!\frac{1}{\rho},\tau_{\mathscr{O}}\!-\!\tau_{\rho}\!>\!\varepsilon\right)\nonumber\\
\label{eq:EstTailProbTauOTauRho} &\leq 3D_{4}\max_{x\in\mathcal{C}}|x|\rho^{-1}+\varepsilon^{-1}\psi(\rho^{-1}).
\end{align}
Next, by (D.7) in~\cite{FlemingSoner:2006}, for some constant $D_{6}>0$ depending on $K$ and $T$,
\begin{align*}
\mathbb{E}^{s,x}\!\left(\left|P_{\tau_{\mathscr{O}}}\right|\left|Z_{\tau_{\mathscr{O}}}\!-Z_{\tau_{\rho}}\right|{\bf 1}_{\{\tau_{\mathscr{O}}-\tau_{\rho}\leq\varepsilon\}}\right)&\leq\left(\mathbb{E}^{s,x}\!\left(\sup_{t\in[s,T]}P_{t}^{2}\right)\right)^{1/2}\left(\mathbb{E}\left(\left|Z_{\tau_{\mathscr{O}}}\!-Z_{\tau_{\rho}}\right|^{2}{\bf 1}_{\{\tau_{\mathscr{O}}-\tau_{\rho}\leq\varepsilon\}}\right)\right)^{1/2}\\
&\leq D_{6}\max_{x\in\mathcal{C}}|x|\left(\mathbb{E}\left(\left|Z_{\tau_{\mathscr{O}}}\!-Z_{\tau_{\rho}}\right|^{2}{\bf 1}_{\{\tau_{\mathscr{O}}-\tau_{\rho}\leq\varepsilon\}}\right)\right)^{1/2}.
\end{align*}
By conditioning on $\mathscr{F}_{\tau_{\rho}}$ and using (D.11) in~\cite{FlemingSoner:2006},
\begin{align*}
\mathbb{E}\left(\left|Z_{\tau_{\mathscr{O}}}\!\!-\!Z_{\tau_{\rho}}\right|^{2}{\bf 1}_{\{\tau_{\mathscr{O}}-\tau_{\rho}\leq\varepsilon\}}\right)&=\int_{s}^{T}\!\int_{\partial\mathscr{O}_{\rho}}\mathbb{E}\left(\left.\left|Z_{\tau_{\mathscr{O}}}\!\!-\!Z_{\tau_{\rho}}\right|^{2}{\bf 1}_{\{\tau_{\mathscr{O}}-\tau_{\rho}\leq\varepsilon\}}\right|(\tau_{\rho},X_{\tau_{\rho}})=(t,y)\right)dF_{\tau_{\rho},X_{\tau_{\rho}}}(t,y)\\
&\leq\int_{s}^{T}\int_{\partial\mathscr{O}_{\rho}}\mathbb{E}^{t,y}\left(\max_{u\in[t,t+\varepsilon]}\left|Z_{u}-y\right|^{2}\right)dF_{\tau_{\rho},X_{\tau_{\rho}}}(t,y)\leq D_{7}\varepsilon,
\end{align*}
where $D_{7}>0$ is a constant depending only on $K$ and $T$. Hence, for any $\varepsilon>0$, there exists a constant $\widetilde{D}>0$ depending on $K$, $T$, $H$, $A$, $\mathscr{N}$, $\mathscr{L}$, $\mathscr{C}$, $\varrho$, $D_{4}$, $D_{5}$, $D_{6}$ and $D_{7}$, such that
\begin{align}\label{eq:EstDiffJRhoJ2}
\mathbb{E}^{s,x}\!\left(\left|P_{\tau_{\mathscr{O}}}\right|\left|Z_{\tau_{\mathscr{O}}}-Z_{\tau_{\rho}}\right|\right)\leq\widetilde{D}\max_{x\in\mathcal{C}}|x|^{3}\left(\varepsilon+\rho^{-1}+\varepsilon^{-1}\psi(\rho^{-1})\right).
\end{align}
Similarly, for any $\varepsilon>0$,
\begin{align}\label{eq:EstDiffJRhoJ3}
\mathbb{E}^{s,x}\!\left(\left|Z_{\tau_{\mathscr{O}}}\right|\left|P_{\tau_{\mathscr{O}}}-P_{\tau_{\rho}}\right|\right)\leq\widetilde{D}\max_{x\in\mathcal{C}}|x|^{3}\left(\varepsilon+\rho^{-1}+\varepsilon^{-1}\psi(\rho^{-1})\right).
\end{align}

Combining \eqref{eq:DecompJRhoJ}, \eqref{eq:EstDiffJRhoJ1}, \eqref{eq:EstDiffJRhoJ2} and \eqref{eq:EstDiffJRhoJ3}, we have shown that, as $\rho\rightarrow\infty$, $J_{\rho,\nu}(s,x;U)$ converges to $J_{\nu}(s,x;U)$, uniformly for all $s\in[0,T]$, $x\in\mathcal{C}$, $\nu\in\mathscr{P}_{[s,T]}$ and $U\in\mathcal{U}_{\nu}[s,T]$. Hence, as $\rho\rightarrow\infty$, $V_{\rho,\nu}(s,x)$ converges to $V_{\nu}(s,y)$, uniformly for all $s\in[0,T]$, $x\in\mathcal{C}$ and $\nu\in\mathscr{P}_{[s,T]}$, which implies that $V_{\nu}\in C([0,T]\times\mathcal{C})$. Since $V_{\rho,\nu}(s,x)=V_{\rho}(s,x)$ for any $(s,x)\in[0,T]\times\mathcal{C}$ and any $\nu\in\mathscr{P}_{[s,T]}$, we have $V_{\nu}(s,x)=V(s,x)$ for all $(s,x)\in[0,T]\times\mathcal{C}$ and $\nu\in\mathscr{P}_{[s,T]}$, and in particular, $V\in C([0,T]\times\mathcal{C})$. Since $\mathcal{C}$ is an arbitrarily chosen compact set, it follows that $V_{\nu}\in C([0,T]\times(\overline{\mathscr{O}}\setminus\{z=0\}))$, that $V\in C([0,T]\times(\overline{\mathscr{O}}\setminus\{z=0\}))$, and that $V_{\nu}(s,x)=V(s,x)$ for all $(s,x)\in[0,T]\times(\overline{\mathscr{O}}\setminus\{z=0\})$ and $\nu\in\mathscr{P}_{[s,T]}$. The continuity of $J_{\nu}(s,x;U)$ on $[0,T]\times\{z=0\}$ (and hence of $V_{\nu}(s,x)$) then follows from \eqref{eq:EstCostFunt} and the fact that $J(s,x;U)=0$ on $[0,T]\times\{z=0\}$. Therefore, $V_{\nu}\in C(\overline{Q}_{T})$, $V\in C(\overline{Q}_{T})$, and $V_{\nu}(s,x)=V(s,x)$ for all $(s,x)\in\overline{Q}_{T}$ and $\nu\in\mathscr{P}_{[s,T]}$. The proof is now complete.\hfill $\Box$

\subsection{The Dynamic Programming Principle}\label{subsec:DynProgPrin}

In order to prove that the value function is a viscosity solution of the HJB equation, and besides the joint continuity, we also need to show that the value function satisfies the so-called Dynamic Programming Principle (cf. (7.2) in~\cite[Section III.7]{FlemingSoner:2006}).
\begin{definition}\label{def:DynProgPrin}
The value function $V$ is said to satisfy the \emph{Dynamic Programming Principle} if, for any $(s,x)\in\overline{Q}_{T}$ and any $(\mathscr{F}_t)_{t\in[s,T]}$-stopping time $\tau$,
\begin{align*}
V(s,x)=\inf_{\substack{U\in\mathcal{U}_{\nu}[s,T] \\ \nu\in\mathscr{P}_{[s,T]}}}\mathbb{E}^{s,x}\left(\int_{s}^{\tau_{\mathscr{O}}\wedge\tau}\mathcal{L}\left(X_{t},U_{t}\right)dt+V\left(\tau_{\mathscr{O}}\wedge\tau,X_{\tau_{\mathscr{O}}\wedge\tau}\right)\right).
\end{align*}
\end{definition}
Indeed, we will verify a stronger version of the traditional dynamic programming principle in this section (cf.~\cite[Definition IV.7.1]{FlemingSoner:2006}).
\begin{definition}\label{def:PropDP}
The value function $V$ is said to satisfy the \emph{property (DP)} if, for any $(s,x)\in\overline{Q}_{T}$,
\begin{itemize}
\item [(i)] for any $\nu\in\mathscr{P}_{[s,T]}$, $U\in\mathcal{U}_{\nu}[s,T]$, and $(\mathscr{F}_t)_{t\in[s,T]}$-stopping time $\tau$,
    \begin{align*}
    V(s,x)\leq\mathbb{E}^{s,x}\left(\int_{s}^{\tau_{\mathscr{O}}\wedge\tau}\mathcal{L}\left(X_{t},U_{t}\right)dt+V\left(\tau_{\mathscr{O}}\wedge\tau,X_{\tau_{\mathscr{O}}\wedge\tau}\right)\right);
    \end{align*}
\item [(ii)] for any $\varepsilon>0$, there exist $\hat{\nu}\in\mathscr{P}_{[s,T]}$ and $\widehat{U}\in\mathcal{U}_{\hat{\nu}}[s,T]$, so that for any $(\mathscr{F}_{t})_{t\in[s,T]}$-stopping time $\tau$,
    \begin{align*}
    V(s,x)+\varepsilon\geq\mathbb{E}^{s,x}\left(\int_{s}^{\tau_{\mathscr{O}}\wedge\tau}\mathcal{L}\left(X_{t},\widehat{U}_{t}\right)dt+V\left(\tau_{\mathscr{O}}\wedge\tau,X_{\tau_{\mathscr{O}}\wedge\tau}\right)\right).
    \end{align*}
\end{itemize}
\end{definition}
Clearly the property (DP) implies the validity of the Dynamic Programming Principle. In the next theorem, we establish the validity of the property (DP) for our value function \eqref{eq:WeakValFunt2}, using the same perturbation scheme as in the proof of Theorem \ref{thm:JointCont}.
\begin{theorem}\label{thm:PropDP}
Under Assumption \ref{assump:Coefficients}, Assumption \ref{assump:SmallBallCond}, and Assumption \ref{assump:SmallBallCondRhon}, the value function $V$, given as in \eqref{eq:WeakValFunt2}, satisfies the property (DP), and thus satisfies the Dynamic Programming Principle.
\end{theorem}

\noindent
\textbf{Proof: Step 1.} We first consider the stochastic control problem \eqref{eq:RhonCostFunt}-\eqref{eq:RhonValFunt2}. In this case, the property (DP) for $V_{\rho}^{(n)}$ was established in \eqref{eq:DPRhon1} and \eqref{eq:DPRhon2}.

\medskip
\noindent
\textbf{Step 2.} We next consider the stochastic control problem \eqref{eq:RhoCostFunt}-\eqref{eq:RhoValFunt2}. For any $(s,x)\in\overline{Q}_{\rho,T}$, $\mu\in\widetilde{\mathscr{P}}_{[s,T]}$ (recalling that there is a one-to-one correspondence between $\widetilde{\mathscr{P}}_{[s,T]}$ and $\mathscr{P}_{[s,T]}$), $U\in\mathcal{V}_{\mu}[s,T]$, and any $(\mathscr{F}_{t})_{t\in[s,T]}$-stopping time $\tau$,
\begin{align}
&\mathbb{E}^{s,x}\left(\left|V_{\rho}\!\left(\tau_{\rho}\wedge\tau,X_{\tau_{\rho}\wedge\tau}^{(\rho)}\right)-V_{\rho}^{(n)}\!\left(\tau_{\rho,n}\wedge\tau,X_{\tau_{\rho,n}\wedge\tau}^{(\rho,n)}\right)\right|\right)\nonumber\\
\label{eq:DecompVRhoVRhon} &\quad\leq\sup_{(s,x)\in\overline{Q}_{\rho,T}}\left|V_{\rho}(s,x)-V_{\rho}^{(n)}(s,x)\right|+\mathbb{E}^{s,x}\left(\left|V_{\rho}\!\left(\tau_{\rho}\wedge\tau,X_{\tau_{\rho}\wedge\tau}^{(\rho)}\right)-V_{\rho}\!\left(\tau_{\rho,n}\wedge\tau,X_{\tau_{\rho,n}\wedge\tau}^{(\rho,n)}\right)\right|\right).
\end{align}
By combining \eqref{eq:DecompJRhoJRhon}-\eqref{eq:EstDiffJRhoJRhon3}, for any $\varepsilon>0$, there exists $N_{1}\in\mathbb{N}$ and a constant $\widetilde{D}_{1}>0$, depending only on $\rho$, $T$, $K$ and $\epsilon$, so that for any $n\geq N_{1}$,
\begin{align}\label{eq:EstSupVRhoVRhon}
\sup_{(s,x)\in\overline{Q}_{\rho,T}}\left|V_{\rho}(s,x)-V_{\rho}^{(n)}(s,x)\right|\leq\widetilde{D}_{1}\varepsilon.
\end{align}
Next, since $V_{\rho}$ is uniformly continuous on $\overline{Q}_{\rho,T}$, there exists $\delta>0$ so that, for any $(t,y),(t',y')\in\overline{Q}_{\rho,T}$ with $|t-t'|\leq\delta$ and $\|y-y'\|\leq\delta$,
\begin{align*}
\left|V_{\rho}(t,y)-V_{\rho}(t',y')\right|\leq\varepsilon.
\end{align*}
Hence,
\begin{align*}
&\mathbb{E}^{s,x}\!\left(\left|V_{\rho}\!\left(\tau_{\rho}\wedge\tau,X_{\tau_{\rho}\wedge\tau}^{(\rho)}\right)-V_{\rho}\!\left(\tau_{\rho,n}\wedge\tau,X_{\tau_{\rho,n}\wedge\tau}^{(\rho,n)}\right)\right|\right)\\
&\quad\leq 2\!\max_{(s,x)\in\overline{Q}_{\rho,T}}\!\left|V_{\rho}(s,x)\right|\mathbb{P}^{s,x}\!\left(\left|X_{\tau_{\rho}}^{(\rho)}\!-\!X_{\tau_{\rho,n}}^{(\rho,n)}\right|>\frac{\delta}{3}\right)+2\!\max_{(s,x)\in\overline{Q}_{\rho,T}}\!\left|V_{\rho}(s,x)\right|\mathbb{P}^{s,x}\!\left(\left|X_{\tau_{\rho}}^{(\rho)}\!-\!X_{\tau}^{(\rho,n)}\right|>\frac{\delta}{3}\right)\\
&\quad\quad\,+2\!\max_{(s,x)\in\overline{Q}_{\rho,T}}\!\left|V_{\rho}(s,x)\right|\mathbb{P}^{s,x}\!\left(\left|X_{\tau}^{(\rho)}\!-\!X_{\tau_{\rho,n}}^{(\rho,n)}\right|>\frac{\delta}{3}\right)+2\!\max_{(s,x)\in\overline{Q}_{\rho,T}}\!\left|V_{\rho}(s,x)\right|\mathbb{P}^{s,x}\!\left(\left|\tau_{\rho}\!-\!\tau_{\rho,n}\right|>\delta\right)+\varepsilon.
\end{align*}
Using arguments similar to those used in obtaining \eqref{eq:EstDiffJRhoJRhon3}, we can show that there exists $N_{2}>0$ and a constant $\widetilde{D}_{2}>0$, depending only on $\rho$, $T$, $K$ and $\epsilon$, so that for any $n\geq N_{2}$,
\begin{align*}
\mathbb{P}^{s,x}\!\left(\left|X_{\tau_{\rho}}^{(\rho)}-X_{\tau}^{(\rho,n)}\right|>\frac{\delta}{3}\right)+\mathbb{P}^{s,x}\!\left(\left|X_{\tau}^{(\rho)}-X_{\tau_{\rho,n}}^{(\rho,n)}\right|>\frac{\delta}{3}\right)\leq\widetilde{D}_{2}\varepsilon.
\end{align*}
Together with \eqref{eq:EstDiffJRhoJRhon2} and \eqref{eq:EstDiffJRhoJRhon3}, there exists $N_{3}\in\mathbb{N}$ and a constant $\widetilde{D}_{3}>0$, depending only on $\rho$, $T$, $K$ and $\epsilon$, so that for any $n\geq N_{3}$,
\begin{align}\label{eq:EstExpVRhoVRhon}
\mathbb{E}^{s,x}\!\left(\left|V_{\rho}\!\left(\tau_{\rho}\wedge\tau,X_{\tau_{\rho}\wedge\tau}^{(\rho)}\right)-V_{\rho}\!\left(\tau_{\rho,n}\wedge\tau,X_{\tau_{\rho,n}\wedge\tau}^{(\rho,n)}\right)\right|\right)\leq\max_{(s,x)\in\overline{Q}_{\rho,T}}\left|V_{\rho}(s,\!x)\right|\widetilde{D}_{3}\varepsilon.
\end{align}
By \eqref{eq:DecompVRhoVRhon}-\eqref{eq:EstExpVRhoVRhon}, for any $(s,x)\in\overline{Q}_{\rho,T}$, $\mu\in\widetilde{\mathscr{P}}_{[s,T]}$, $U\in\mathcal{V}_{\mu}[s,T]$, any $(\mathscr{F}_{t})_{t\in[s,T]}$-stopping time $\tau$,
\begin{align}\label{eq:LimnSupVRhoVRhon}
\lim_{n\rightarrow\infty}\sup_{\tau}\mathbb{E}^{s,x}\left(\left|V_{\rho}\!\left(\tau_{\rho}\wedge\tau,X_{\tau_{\rho}\wedge\tau}^{(\rho)}\right)-V_{\rho}^{(n)}\!\left(\tau_{\rho,n}\wedge\tau,X_{\tau_{\rho,n}\wedge\tau}^{(\rho,n)}\right)\right|\right)=0.
\end{align}
Moreover, by \eqref{eq:EstDiffJRhoJRhon1} and \eqref{eq:EstDiffJRhoJRhon2}, for any $(s,x)\in\overline{Q}_{\rho,T}$, $\mu\in\widetilde{\mathscr{P}}_{[s,T]}$, $U\in\mathcal{V}_{\mu}[s,T]$ and any $(\mathscr{F}_{t})_{t\in[s,T]}$-stopping time $\tau$, there exist $N_{4}\in\mathbb{N}$ and a constant $\widetilde{D}_{4}>0$, depending only on $\kappa$, $\mathscr{L}$, $\gamma$, $\rho$, $T$, $K$ and $\epsilon$, so that for any $n\geq N_{4}$,
\begin{align}
&\mathbb{E}^{s,x}\left(\left|\int_{s}^{\tau_{\rho}\wedge\tau}\mathcal{L}_{\rho}\!\left(X_{t}^{(\rho)},U_{t}\right)dt-\int_{s}^{\tau_{\rho,n}\wedge\tau}\mathcal{L}_{\rho}\!\left(X_{t}^{(\rho,n)},U_{t}\right)dt\right|\right)\nonumber\\ &\quad\leq\mathbb{E}^{s,x}\left(\int_{s}^{\tau_{\rho}\wedge\tau_{\rho,n}}\left|\mathcal{L}_{\rho}\!\left(X_{t}^{(\rho)},U_{t}\right)-\mathcal{L}_{\rho}\!\left(X_{t}^{(\rho,n)},U_{t}\right)\right|dt\right)+\kappa\mathscr{L}^{\gamma}\rho\,\mathbb{E}^{s,x}\left(\left|\tau_{\rho}-\tau_{\rho,n}\right|\right)\nonumber\\
\label{eq:EstLRhoRhon} &\quad\leq\widetilde{D}_{4}\varepsilon.
\end{align}

Combining \eqref{eq:DPRhon1} with \eqref{eq:LimnSupVRhoVRhon} and \eqref{eq:EstLRhoRhon}, for any $(s,x)\in\overline{Q}_{\rho,T}$, $\mu\in\widetilde{\mathscr{P}}_{[s,T]}$, $U\in\mathcal{V}_{\mu}[s,T]$, any $(\mathscr{F}_{t})_{t\in[s,T]}$-stopping time $\tau$,
\begin{align}\label{eq:DPRho1}
V_{\rho}(s,x)\leq\mathbb{E}^{s,x}\left(\int_{s}^{\tau_{\rho}\wedge\tau}\mathcal{L}_{\rho}\!\left(X_{t}^{(\rho)},U_{t}\right)dt+V_{\rho}\left(\tau_{\rho}\wedge\tau,X_{\tau_{\rho}\wedge\tau}^{(\rho)}\right)\right).
\end{align}
Moreover, by \eqref{eq:LimnSupVRhoVRhon} and \eqref{eq:EstLRhoRhon}, pick $n\in\mathbb{N}$ large enough, so that for any $(s,x)\in\overline{Q}_{\rho,T}$, $\mu\in\widetilde{\mathscr{P}}_{[s,T]}$ and $U\in\mathcal{V}_{\mu}[s,T]$,
\begin{align*}
\sup_{\tau}\mathbb{E}^{s,x}\!\!\left(\left|\int_{s}^{\tau_{\rho}\wedge\tau}\!\!\!\!\mathcal{L}_{\rho}\!\!\left(\!X_{t}^{(\rho)}\!,U_{t}\!\right)\!dt\!+\!V_{\rho}\!\left(\!\tau_{\rho}\!\wedge\!\tau,X_{\tau_{\rho}\wedge\tau}^{(\rho)}\!\right)\!-\!\!\int_{s}^{\tau_{\rho,n}\wedge\tau}\!\!\!\!\mathcal{L}_{\rho}\!\!\left(\!X_{t}^{(\rho,n)}\!,U_{t}\!\right)\!dt\!-\!V_{\rho}^{(n)}\!\left(\!\tau_{\rho,n}\!\wedge\!\tau,X_{\tau_{\rho,n}\wedge\tau}^{(\rho,n)}\!\right)\right|\right)\!\leq\!\frac{\varepsilon}{2}.
\end{align*}
For this choice of $n$, by \eqref{eq:DPRhon2}, there exist $\hat{\mu}\in\widetilde{\mathscr{P}}_{[s,T]}$ and $\widehat{U}\in\mathcal{V}_{\mu}[s,T]$, such that for any $(\mathscr{F}_{t})_{t\in[s,T]}$-stopping times $\tau$,
\begin{align*}
V_{\rho}^{(n)}(s,x)+\frac{\varepsilon}{2}\geq\mathbb{E}^{s,x}\left(\int_{s}^{\tau_{\rho,n}\wedge\tau}\mathcal{L}_{\rho}\left(X_{t}^{(\rho,n)},\widehat{U}_{t}\right)dt+V_{\rho}^{(n)}\left(\tau_{\rho,n}\wedge\tau,X_{\tau_{\rho,n}\wedge\tau}^{(\rho,n)}\right)\right).
\end{align*}
Note that by Remark \ref{rem:EquiStoSym}, we can take $\hat{\nu}\in\mathscr{P}_{[s,T]}$ by omitting the last component $\widetilde{W}$ of the six-tuple $\hat{\mu}$, and hence $\widehat{U}\in\mathcal{U}_{\nu}[s,T]$. Therefore, we find $\hat{\nu}\in\mathscr{P}_{[s,T]}$ and $\widehat{U}\in\mathcal{U}_{\nu}[s,T]$, so that for any $(\mathscr{F}_{t})_{t\in[s,T]}$-stopping times $\tau$,
\begin{align}\label{eq:DPRho2}
V_{\rho}(s,x)+\varepsilon\geq\mathbb{E}^{s,x}\left(\int_{s}^{\tau_{\rho}\wedge\tau}\mathcal{L}_{\rho}\left(X_{t}^{(\rho)},\widehat{U}_{t}\right)dt+V_{\rho}\left(\tau_{\rho}\wedge\tau,X_{\tau_{\rho}\wedge\tau}^{(\rho)}\right)\right).
\end{align}
Therefore, the value function $V_{\rho}$ satisfies the property (DP).

\medskip
\noindent
\textbf{Step 3.} Finally, we consider the stochastic control problem \eqref{eq:WeakCostFunt}-\eqref{eq:WeakValFunt2}, and establish the property (DP) for $V$. We first notice that the property (DP) is trivial when $x=(p,0,\theta)$. Now for any $(s,x)\in[0,T]\times(\overline{Q}_{T}\setminus\{z=0\})$, pick $\rho_{0}>0$ large enough so that $x\in\overline{\mathscr{O}}_{\rho_{0}}$ (and hence $|x|\leq\rho_{0}$). For any $\nu\in\mathscr{P}_{[s,T]}$, $U\in\mathcal{U}_{\nu}[s,T]$, any $(\mathscr{F}_{t})_{t\in[s,T]}$-stopping time $\tau$, and any $\rho\geq\rho_{0}$,
\begin{align}
&\mathbb{E}^{s,x}\left(\left|V_{\rho}\left(\tau_{\rho}\wedge\tau,X_{\tau_{\rho}\wedge\tau}^{(\rho)}\right)-V\left(\tau_{\mathscr{O}}\wedge\tau,X_{\tau_{\mathscr{O}}\wedge\tau}\right)\right|\right)\nonumber\\
&\quad\leq\mathbb{E}^{s,x}\left(\left|V_{\rho}\left(\tau_{\rho}\wedge\tau,X_{\tau_{\rho}\wedge\tau}^{(\rho)}\right)-V\left(\tau_{\mathscr{O}}\wedge\tau,X_{\tau_{\mathscr{O}}\wedge\tau}\right)\right|{\bf 1}_{\{\sup_{t\in[s,T]}|X_{t}|>\rho_{0}^{2}\}}\right)\nonumber\\
\label{eq:DecompVRhoV} &\quad\quad\,+\mathbb{E}^{s,x}\left(\left|V_{\rho}\left(\tau_{\rho}\wedge\tau,X_{\tau_{\rho}\wedge\tau}^{(\rho)}\right)-V\left(\tau_{\mathscr{O}}\wedge\tau,X_{\tau_{\mathscr{O}}\wedge\tau}\right)\right|{\bf 1}_{\{\sup_{t\in[s,T]}|X_{t}|\leq\rho_{0}^{2}\}}\right).
\end{align}
To estimate the first expectation in \eqref{eq:DecompVRhoV}, by Lemma \ref{lem:WellDefStoContProb}, there exists a constant $\widetilde{K}>0$, depending only on $\kappa$, $\mathscr{L}$, $\gamma$, $T$, $K$, $H$, $A$, $\mathscr{N}$, $\mathscr{C}$, $\alpha$, $\beta$ and $\varrho$, so that for any $(s,x)\in\overline{Q}_{T}$,
\begin{align}\label{eq:EstV}
\left|V(s,x)\right|\leq\widetilde{K}\left(1+z+z^{2}+p^{2}\right),
\end{align}
where $x=(p,z,\theta)$. The same estimate holds for $V_{\rho}$. Hence, by (D.7) in~\cite{FlemingSoner:2006},
\begin{align}
&\mathbb{E}^{s,x}\left(\left|V_{\rho}\left(\tau_{\rho}\wedge\tau,X_{\tau_{\rho}\wedge\tau}^{(\rho)}\right)-V\left(\tau_{\mathscr{O}}\wedge\tau,X_{\tau_{\mathscr{O}}\wedge\tau}\right)\right|{\bf 1}_{\{\sup_{t\in[s,T]}|X_{t}|>\rho_{0}^{2}\}}\right)\nonumber\\
&\quad\leq 2\widetilde{K}\,\mathbb{E}^{s,x}\left(\left(1+\sup_{t\in[s,T]}\left|X_{t}\right|+\sup_{t\in[s,T]}\left|X_{t}\right|^{2}\right){\bf 1}_{\{\sup_{t\in[s,T]}|X_{t}|>\rho_{0}^{2}\}}\right)\nonumber\\
&\quad\leq 2\widetilde{K}\,\mathbb{P}^{s,x}\!\left(\sup_{t\in[s,T]}|X_{t}|>\rho_{0}^{2}\right)+2\widetilde{K}\,\sqrt{\mathbb{E}^{s,x}\!\left(\sup_{t\in[s,T]}\left|X_{t}\right|^{2}\right)\mathbb{P}^{s,x}\!\left(\sup_{t\in[s,T]}|X_{t}|>\rho_{0}^{2}\right)}\nonumber\\
&\quad\quad\,+2\widetilde{K}\,\sqrt{\mathbb{E}^{s,x}\!\left(\sup_{t\in[s,T]}\left|X_{t}\right|^{4}\right)\mathbb{P}^{s,x}\!\left(\sup_{t\in[s,T]}|X_{t}|>\rho_{0}^{2}\right)}\nonumber\\ &\quad\leq\widetilde{D}_{5}\rho_{0}^{-2}\left(1+|x|\right)+\widetilde{D}_{5}\sqrt{\rho_{0}^{-6}\left(1+|x|^{2}\right)\left(1+|x|^{3}\right)}+\widetilde{D}_{5}\sqrt{\rho_{0}^{-10}\left(1+|x|^{4}\right)\left(1+|x|^{5}\right)}\nonumber\\
\label{eq:EstDifVRhoV1} &\quad\leq\frac{6\widetilde{D}_{5}}{\sqrt{\rho_{0}}},
\end{align}
where $\widetilde{D}_{5}$ is a constant depending only on $\kappa$, $\mathscr{L}$, $\gamma$, $T$, $K$, $H$, $A$, $\mathscr{N}$, $\mathscr{C}$, $\alpha$, $\beta$ and $\varrho$. Note that since $V$ is uniformly continuous on $[0,T]\times(\overline{Q}_{T}\cap\{|x|\leq\rho_{0}\})$, there exists $\delta\in(0,\varepsilon)$, such that for any $(t,y),(t',y')\in [0,T]\times(\overline{Q}_{T}\cap\{|x|\leq\rho_{0}\})$ with $|t-t'|\leq\delta$ and $|y-y'|\leq\delta$,
\begin{align*}
\left|V(t,y)-V(t',y')\right|\leq\varepsilon.
\end{align*}
Hence, together with \eqref{eq:EstV}, and noting that $X$ and $X^{(\rho)}$ are identical up to $\tau_{\rho}$,
\begin{align*}
&\mathbb{E}^{s,x}\left(\left|V_{\rho}\left(\tau_{\rho}\wedge\tau,X_{\tau_{\rho}\wedge\tau}^{(\rho)}\right)-V\left(\tau_{\mathscr{O}}\wedge\tau,X_{\tau_{\mathscr{O}}\wedge\tau}\right)\right|{\bf 1}_{\{\sup_{t\in[s,T]}|X_{t}|\leq\rho_{0}^{2}\}}\right)\\
&\quad\leq\varepsilon+6\widetilde{K}\rho_{0}^{2}\left(\mathbb{P}^{s,x}\left(\tau_{\mathscr{O}}-\tau_{\rho}>\delta\right)+\mathbb{P}^{s,x}\left(\left|X_{\tau_{\rho}\wedge\tau}-X_{\tau_{\mathscr{O}}\wedge\tau}\right|>\delta,\,\sup_{t\in[s,T]}|X_{t}|\leq\rho_{0}^{2}\right)\right).
\end{align*}
Above, the first probability is already estimated in \eqref{eq:EstTailProbTauOTauRho} (with $\mathcal{C}=\overline{\mathscr{O}}_{\rho_{0}}$). Moreover, by \eqref{eq:EstTailProbTauOTauRho} and (D.12) in~\cite{FlemingSoner:2006}, for some $\widetilde{D}_{6}>0$, depending only on $\mathscr{L}$, $T$, $K$, $H$, $A$, $\mathscr{N}$, $\mathscr{C}$, $\alpha$, $\beta$ and $\varrho$,
\begin{align*}
&\mathbb{P}^{s,x}\bigg(\left|X_{\tau_{\rho}\wedge\tau}-X_{\tau_{\mathscr{O}}\wedge\tau}\right|>\delta,\,\sup_{t\in[s,T]}|X_{t}|\leq\rho_{0}^{2}\bigg)\\
&\quad\leq\mathbb{P}^{s,x}\bigg(\left|X_{\tau_{\rho}\wedge\tau}-X_{\tau_{\mathscr{O}}\wedge\tau}\right|>\delta,\,\tau_{\mathscr{O}}-\tau_{\rho}\leq\delta^{3},\,\sup_{t\in[s,T]}|X_{t}|\leq\rho_{0}^{2}\bigg)+\mathbb{P}^{s,x}\!\left(\tau_{\mathscr{O}}-\tau_{\rho}>\delta^{3}\right)\\
&\quad\leq\widetilde{D}_{6}\,\rho_{0}^{2}\,\delta+\widetilde{D}_{6}\rho_{0}\rho^{-1}+\delta^{-3}\psi(\rho^{-1}).
\end{align*}
Hence, the second expectation in \eqref{eq:DecompVRhoV} is bounded by
\begin{align}
&\mathbb{E}^{s,x}\left(\left|V_{\rho}\left(\tau_{\rho}\wedge\tau,X_{\tau_{\rho}\wedge\tau}^{(\rho)}\right)-V\left(\tau_{\mathscr{O}}\wedge\tau,X_{\tau_{\mathscr{O}}\wedge\tau}\right)\right|{\bf 1}_{\{\sup_{t\in[s,T]}|X_{t}|\leq\rho_{0}^{2}\}}\right)\nonumber\\
\label{eq:EstDiffVRhoV2} &\quad\leq\varepsilon+6\widetilde{K}\rho_{0}^{2}\left(\widetilde{D}_{6}\rho_{0}^{2}\delta+\widetilde{D}_{6}\rho_{0}\rho^{-1}+\delta^{-3}\psi(\rho^{-1})+3D_{4}\rho_{0}\rho^{-1}+\delta^{-1}\psi(\rho^{-1})\right).
\end{align}
Moreover, by \eqref{eq:DecompJRhoJ} and \eqref{eq:EstDiffJRhoJ1} (with $\mathcal{C}=\overline{\mathscr{O}}_{\rho_{0}}$), there exists a constant $\widetilde{D}_{7}>0$, depending on $\mathscr{L}$, $T$, $K$, $H$, $A$, $\mathscr{N}$, $\mathscr{C}$, $\alpha$, $\beta$ and $\varrho$, such that for any $\nu\in\mathscr{P}_{[s,T]}$, any $U\in\mathcal{U}_{\nu}[s,T]$, and any $(\mathscr{F}_{t})_{t\in[s,T]}$-stopping time $\tau$,
\begin{align}
\mathbb{E}^{s,x}\left(\left|\int_{s}^{\tau_{\rho}\wedge\tau}\mathcal{L}_{\rho}\!\left(X_{t}^{(\rho)},U_{t}\right)dt-\int_{s}^{\tau_{\mathscr{O}}\wedge\tau}\mathcal{L}\left(X_{t},U_{t}\right)dt\right|\right)&\leq\mathbb{E}^{s,x}\left(\int_{\tau_{\rho}\wedge\tau}^{\tau_{\mathscr{O}}\wedge\tau}\mathcal{L}\left(X_{t},U_{t}\right)dt\right)\nonumber\\
\label{eq:EstDiffLRhoL} &\leq\widetilde{D}_{7}\,\rho_{0}^{3}\sqrt{\rho^{-1}+\delta^{-1}\psi(\rho^{-1})+\delta}.
\end{align}

Combining \eqref{eq:DecompVRhoV}-\eqref{eq:EstDiffLRhoL}, for any $(s,x)\in[0,T]\times(\overline{Q}_{T}\setminus\{z=0\})$ and any $\varepsilon>0$, choose first $\rho_{0}>0$ large enough so that $x\in\overline{\mathscr{O}}_{\rho_{0}}$ and that $\sqrt{\rho_{0}}\geq\varepsilon^{-1}$. Next, choose $\delta>0$ small enough, so that $\delta\leq\varepsilon^{2}/\rho_{0}^{6}$. Finally, choose $\rho_{1}>\rho_{0}$, large enough, so that
\begin{align*}
\rho_{1}\geq\frac{\rho_{0}^{6}}{\varepsilon^{2}},\quad\frac{\rho_{0}^{2}}{\delta^{3}}\psi(\rho_{1}^{-1})\leq\varepsilon,\quad\frac{\rho_{0}^{6}}{\delta}\psi(\rho_{1}^{-1})\leq\varepsilon^{2}.
\end{align*}
Then, for any $\rho\geq\rho_{1}$, for any $\nu\in\mathscr{P}_{[s,T]}$ and any $U\in\mathcal{U}_{\nu}[s,T]$, there exists a constant $\widetilde{D}>0$ depending only on $\mathscr{L}$, $T$, $K$, $H$, $A$, $\mathscr{N}$, $\mathscr{C}$, $\alpha$, $\beta$ and $\varrho$, so that
\begin{align}\label{eq:EstSupVRhoVLRhoL}
\sup_{\tau}\mathbb{E}^{s,x}\!\!\left(\!\left|V_{\rho}\!\!\left(\!\tau_{\rho}\!\wedge\!\tau,X_{\tau_{\rho}\wedge\tau}^{(\rho)}\!\right)\!-\!V\!\left(\!\tau_{\mathscr{O}}\!\wedge\!\tau,X_{\tau_{\mathscr{O}}\wedge\tau}\!\right)\!\right|\!+\!\left|\int_{s}^{\tau_{\rho}\wedge\tau}\!\!\!\!\!\mathcal{L}_{\rho}\!\!\left(\!X_{t}^{(\rho)}\!,U_{t}\!\right)\!dt\!-\!\!\!\int_{s}^{\tau_{\mathscr{O}}\wedge\tau}\!\!\!\!\!\mathcal{L}\!\left(X_{t},\!U_{t}\right)\!dt\right|\right)\!\!\leq\!\widetilde{D}\varepsilon.
\end{align}
The validity of the property (DP) for $V$ then follows immediately from \eqref{eq:DPRho1}, \eqref{eq:DPRho2} and \eqref{eq:EstSupVRhoVLRhoL}.\hfill $\Box$

\subsection{Existence of Viscosity Solutions}\label{subsec:ExistVisSols}

Now that we have established the joint continuity and the Dynamic Programming Principle for the value function $V$, it is time to show that the it is indeed a viscosity solution of the HJB equation \eqref{eq:HJB} with terminal/boundary condition \eqref{eq:HJBTermBound}.
\begin{theorem}\label{thm:ExistVisSols}
Under Assumption \ref{assump:Coefficients}, Assumption \ref{assump:SmallBallCond}, and Assumption \ref{assump:SmallBallCondRhon}, the value function $V$, as given in \eqref{eq:WeakValFunt2}, is a viscosity solution of the HJB equation \eqref{eq:HJB} with terminal/boundary condition
\eqref{eq:HJBTermBound}.
\end{theorem}

\noindent
\textbf{Proof:} The proof is very similar to that of~\cite[Theorem 4.5.2]{YongZhou:1999}, and here we only present the outline. The boundary/terminal condition is clearly satisfied. For any $\varphi\in C^{1,2}(\overline{Q}_{T})$, let $V-\varphi$ attain a local maximum at some $(\bar{s},\bar{x})\in Q_{T}$. Without loss of generality, we can assume $\varphi\in C^{1,2}_{b}(\overline{Q}_{T})$, i.e., all derivatives of $\varphi$ are bounded in $\overline{Q}_{T}$. Fix any $u\in\mathscr{U}$, and consider the constant control $U_{t}\equiv u$, $t\in[\bar{s},T]$. By the property (DP)-(i), It\^{o}'s formula and the dominated convergence theorem, for any reference stochastic system $\nu\in\mathscr{P}_{[\bar{s},T]}$, any $s>\bar{s}$ with $s-\bar{s}>0$ small enough,
\begin{align*}
0&\leq\frac{1}{s-\bar{s}}\mathbb{E}^{\bar{s},\bar{x}}\left(V(\bar{s},\bar{x})-\varphi(\bar{s},\bar{x})-V\left(s\wedge\tau_{\mathscr{O}},X_{s\wedge\tau_{\mathscr{O}}}\right)+\varphi\left(s\wedge\tau_{\mathscr{O}},X_{s\wedge\tau_{\mathscr{O}}}\right)\right)\\
&\leq\frac{1}{s-\bar{s}}\mathbb{E}^{\bar{s},\bar{x}}\left(\int_{\bar{s}}^{s\wedge\tau_{\mathscr{O}}}\mathcal{L}\left(X_{t},u\right)dt-\varphi(\bar{s},\bar{x})+\varphi\left(s\wedge\tau_{\mathscr{O}},X_{s\wedge\tau_{\mathscr{O}}}\right)\right)\\
&\quad\longrightarrow\mathcal{L}(\bar{x},u)+\varphi_{s}(\bar{s},\bar{x})+\vec{f}(\bar{s},\bar{x},u)\cdot D_{x}\varphi(\bar{s},\bar{x})+\frac{1}{2}\text{tr}\left(\vec{a}(\bar{s},\bar{x},u)D_{xx}\varphi(\bar{s},\bar{x})\right),\quad s\downarrow\bar{s}.
\end{align*}
Hence,
\begin{align*}
-\varphi_{s}(\bar{s},\bar{x})+\mathcal{H}\left(\bar{s},\bar{x},u,D_{x}\varphi(\bar{s},\bar{x}),D_{xx}\varphi(\bar{s},\bar{x})\right)\leq 0.
\end{align*}
On the other hand, let $V-\varphi$ attain a local minimum at some $(\bar{s},\bar{x})\in Q_{T}$. By the property (DP)-(ii), for any $\varepsilon>0$,  and $s>\bar{s}$ with $s-\bar{s}>0$ small enough, there exists $\hat{\nu}\in\mathscr{P}_{[s,T]}$ and $\widehat{U}\in\mathcal{U}_{\hat{\nu}}[\bar{s},T]$, such that
\begin{align*}
0&\geq\mathbb{E}^{\bar{s},\bar{x}}\left(V(\bar{s},\bar{x})-\varphi(\bar{s},\bar{x})-V\left(s\wedge\tau_{\mathscr{O}},X_{s\wedge\tau_{\mathscr{O}}}\right)+\varphi\left(s\wedge\tau_{\mathscr{O}},X_{s\wedge\tau_{\mathscr{O}}}\right)\right)\\
&\geq
-\varepsilon(s-\bar{s})+\mathbb{E}_{\bar{s},\bar{x}}\left(\int_{\bar{s}}^{s\wedge\tau_{\mathscr{O}}}\mathcal{L}\left(X_{t},\widehat{U}_{t}\right)dt+\varphi\left(s\wedge\tau_{\mathscr{O}},X_{s\wedge\tau_{\mathscr{O}}}\right)-\varphi(\bar{s},\bar{x})\right).
\end{align*}
Hence, by It\^{o}'s formula and the dominated convergence theorem,
\begin{align*}
-\varepsilon&\leq\frac{1}{s-\bar{s}}\mathbb{E}^{\bar{s},\bar{x}}\left(\int_{\bar{s}}^{s\wedge\tau_{\mathscr{O}}}\!\left(-\varphi_{s}(t,X_{t})\!-\!\vec{f}\!\left(t,X_{t},\widehat{U}_{t}\right)\!\cdot\! D_{x}\varphi(t,X_{t})\!-\!\frac{1}{2}\text{tr}\left(\vec{a}\left(t,X_{t},\widehat{U}_{t}\right)\!D_{xx}\varphi(t,X_{t})\right)\!\right)\!dt\!\right)\\
&\leq\frac{1}{s-\bar{s}}\mathbb{E}^{\bar{s},\bar{x}}\left(\int_{\bar{s}}^{s\wedge\tau_{\mathscr{O}}}\left(-\varphi_{s}(t,X_{t})+\mathcal{H}\left(t,X_{t},D_{x}\varphi(t,X_{t}),D_{xx}\varphi(t,X_{t})\right)\right)dt\right)\\
&\quad\longrightarrow-\varphi_{s}(\bar{s},\bar{x})+\mathcal{H}\left(\bar{s},\bar{x},D_{x}\varphi(\bar{s},\bar{x}),D_{xx}\varphi(\bar{s},\bar{x})\right),\quad s\downarrow\bar{s}.
\end{align*}
Since $\varepsilon>0$ is arbitrary,
\begin{align*}
-\varphi_{s}(\bar{s},\bar{x})+\mathcal{H}\left(\bar{s},\bar{x},D_{x}\varphi(\bar{s},\bar{x}),D_{xx}\varphi(\bar{s},\bar{x})\right)\geq 0,
\end{align*}
which completes the proof.\hfill $\Box$

\section{Uniqueness of the Viscosity Solution}\label{sec:Uniqueness}

In this section, we establish the comparison principle for viscosity subsolutions and supersolutions to \eqref{eq:HJB} with terminal/boundary condition \eqref{eq:HJBTermBound}. This, together with Theorem \ref{thm:ExistVisSols}, shows that the value function \eqref{eq:WeakValFunt2} is the unique viscosity solution with polynomial growth (recalling Lemma \ref{lem:WellDefStoContProb}) to \eqref{eq:HJB} with terminal/boundary condition \eqref{eq:HJBTermBound}.
\begin{theorem}\label{thm:ComparisonPrin}
Let Assumption \ref{assump:Coefficients} be satisfied. Let $\mathcal{W}_{1}$ and $\mathcal{W}_{2}$ be, respectively, any subsolution and supersolution to \eqref{eq:HJB}, both of which satisfies the equality of the boundary/terminal condition \eqref{eq:HJBTermBound}. Moreover, let $\mathcal{W}_{1}$ and $\mathcal{W}_{2}$ satisfy a polynomial growth condition in the space variable, i.e., for any $(s,x)\in\overline{Q}_{T}$,
\begin{align}\label{eq:PolyAssumpWi}
\left|\mathcal{W}_{i}(s,x)\right|\leq K_{0}\left(1+|x|^{m}\right),\quad i=1,2,
\end{align}
for some constant $K_{0}>0$ and $m\in\mathbb{N}$. Then, $\mathcal{W}_{1}(s,x)\leq\mathcal{W}_{2}(s,x)$, for any $(s,x)\in\overline{Q}_{T}$. In particular, the value function $V$, defined in \eqref{eq:WeakValFunt2}, is the unique viscosity solution to \eqref{eq:HJB} with terminal/boundary condition \eqref{eq:HJBTermBound}, having polynomial growth in the space variable.
\end{theorem}
\begin{remark}\label{rem:TermBoundCondEqu}
Above, both the subsolution and the supersolution are assumed to satisfy the boundary/terminal condition with equality. This avoids to appeal to extra conditions such as uniform continuity on boundary/terminal values of the solutions.
\end{remark}
The proof of Theorem \ref{thm:ComparisonPrin} relies mainly on the following remarkable result known as \emph{Ishii's lemma}. To state the result, we first introduce some more notations. For any locally compact subset $\mathcal{Q}\subseteq\mathbb{R}^{4}$, let $\text{USC}(\mathcal{Q})$ (respectively, $\text{LSC}(\mathcal{Q})$) be the collection of all real-value upper (respectively, lower) semicontinuous functions on $\mathcal{Q}$. For $x\in\mathcal{Q}$, and $v\in\text{USC}(\mathcal{Q})$, let
\begin{align*}
\mathcal{J}^{2}_{+}v(x)&:=\left\{(r,G)\in\mathbb{R}^{4}\times\mathscr{S}^{4}:\,\,v(x+h)-v(x)\leq r\cdot h+\!\frac{1}{2}h\cdot Gh+o(|h|^{2})\right\}\\
&\,\,=\left\{\left(D_{x}\phi(x),D_{xx}\phi(x)\right):\,\,\phi\in C^{2}(\mathcal{Q}),\,v-\phi\,\,\,\text{has a local maximum at}\,\,x\right\},
\end{align*}
and
\begin{align*}
\overline{\mathcal{J}}^{2}_{+}v(x):=\left\{(r,G)\!\in\!\mathbb{R}^{4}\!\times\!\mathscr{S}^{4}:\,\exists x_{n}\in\mathcal{Q},\,(r_{n},\!G_{n})\in\mathcal{J}^{2,+}v(x_{n}),\,\lim_{n\rightarrow\infty}x_{n}\!=\!x,\,\lim_{n\rightarrow\infty}(r_{n},\!G_{n})\!=\!(r,G)\right\}.
\end{align*}
Similarly, for $x\in\mathcal{Q}$ and $v\in\text{LSC}(\mathcal{Q})$, let
\begin{align*}
\mathcal{J}^{2}_{-}v(x)&:=\left\{(r,G)\in\mathbb{R}^{4}\times\mathscr{S}^{4}:\,\,v(x+h)-v(x)\geq r\cdot h+\frac{1}{2}h\cdot Gh+o(|h|^{2})\right\}\\
&\,\,=\left\{\left(D_{x}\phi(x),D_{xx}\phi(x)\right):\,\,\phi\in C^{2}(\mathcal{Q}),\,v-\phi\,\,\,\text{has a local minimum at}\,\,x\right\},
\end{align*}
and
\begin{align*}
\overline{\mathcal{J}}^{2}_{-}v(x):=\left\{(r,G)\!\in\!\mathbb{R}^{4}\!\times\!\mathscr{S}^{4}:\,\exists x_{n}\in\mathcal{Q},\,(r_{n},\!G_{n})\in\mathcal{J}^{2,-}v(x_{n}),\,\lim_{n\rightarrow\infty}x_{n}\!=\!x,\,\lim_{n\rightarrow\infty}(r_{n},\!G_{n})\!=\!(r,G)\right\}.
\end{align*}

The following version of Ishii's Lemma is taken from~\cite[Theorem 3.2]{CrandallIshiiLions:1992}.
\begin{theorem}\label{thm:IshiiLemma}
For $i=1,\ldots,k$, let $\mathcal{Q}_{i}$ be a locally compact subsets of $\mathbb{R}^{N_{i}}$, where $N_{i}\in\mathbb{N}$. Let $v_{i}\in\text{USC}(\mathcal{Q}_{i})$, and $\varphi\in C^{2}(\mathcal{Q})$, where $\mathcal{Q}:=\mathcal{Q}_{1}\times\cdots\times\mathcal{Q}_{k}$. For $x=(x_{1},\ldots,x_{k})\in\mathcal{Q}$, let
\begin{align*}
v(x):=v_{1}(x_{1})+\cdots+v_{k}(x_{k}).
\end{align*}
Assume that $v-\varphi$ achieves a local maximum at $\bar{x}=(\bar{x}_{1},\ldots,\bar{x}_{k})\in\mathcal{Q}$ then, for any $\varepsilon>0$, there exists $G_{i}\in\mathscr{S}^{N_{i}}$, such that $(D_{x_{i}}\varphi(\bar{x}),G_{i})\in\overline{\mathcal{J}}^{2}_{+}v_{i}(\bar{x}_{i})$, for each $i=1,\ldots,k$, and such that the block diagonal matrix with entries $G_{i}$, $i=1,\ldots,k$, satisfies
\begin{align*}
\left(\begin{array}{ccc} G_{1} & \cdots & 0 \\ \vdots & \ddots & \vdots \\ 0 & \cdots & G_{k} \end{array}\right)\leq G_{\varphi}+\varepsilon\,G_{\varphi}^{2},
\end{align*}
where $G_{\varphi}:=D_{xx}\varphi(\bar{x})\in\mathscr{S}_{+}^{N}$, and where $N:=N_{1}+\cdots+N_{k}$.
\end{theorem}

\noindent
\textbf{Proof of Theorem \ref{thm:ComparisonPrin}:} We argue by contradiction. Suppose that there exist $(s_{0},x_{0})\in\overline{Q}_{T}$ and $\gamma_{0}>0$ so that
\begin{align*}
\mathcal{W}_{1}(s_{0},x_{0})-\mathcal{W}_{2}(s_{0},x_{0})\geq\gamma_{0}>0.
\end{align*}
Choose $r>0$ and $\delta>0$ small enough so that
\begin{align}\label{eq:MaxPsiGamma0}
\mathcal{W}_{1}(s_{0},x_{0})-\mathcal{W}_{2}(s_{0},x_{0})-\frac{2r}{t_{0}}-2\delta\,e^{-T}\left(|x_{0}|^{m+1}+1\right)>\frac{\gamma_{0}}{2}.
\end{align}
Here, without loss of generality, we assume that $s_{0}>0$. Otherwise we can replace $r/t_{0}$ by $r/(T-s_{0})$ and the argument is similar.

\medskip
\noindent
\textbf{Step 1.} For $\varepsilon_{1}>0$, $\varepsilon_{2}>0$ and $\rho>1$, let
\begin{align*}
\Psi(s,x,t,y):=\mathcal{W}_{1}(s,x)-\mathcal{W}_{2}(t,y)-\frac{\delta}{e^{\rho s}}\!\left(|x|^{m+1}\!+\!1\right)-\frac{\delta}{e^{\rho t}}\!\left(|y|^{m+1}\!+\!1\right)-\frac{|x\!-\!y|^{2}}{2\varepsilon_{1}}-\frac{|s\!-\!t|^{2}}{2\varepsilon_{2}}-\frac{r}{s}-\frac{r}{t},
\end{align*}
for $(s,x,t,y)\in\overline{Q}_{T}^{2}$. We claim that $\Psi$ attains its maximum in the interior of $Q_{T}^{2}$. To see this, let
\begin{align*}
&m_{0}:=\lim_{\eta\rightarrow 0}\,\lim_{\varepsilon\rightarrow
0}\sup_{\substack{|t-s|<\varepsilon \\ |x-y|<\eta}}\left(\mathcal{W}_{1}(s,x)-\mathcal{W}_{2}(t,y)-\frac{r}{s}-\frac{r}{t}-\frac{\delta}{e^{\rho s}}\left(|x|^{m+1}+1\right)-\frac{\delta}{e^{\rho t}}\left(|y|^{m+1}+1\right)\right),\\
&m_{1}(\varepsilon_{1},\varepsilon_{2}):=\sup_{(s,x,t,y)\in\overline{Q}_{T}^{2}}\Psi(s,x,t,y),\\
&m_{2}(\varepsilon_{2}):=\lim_{\eta\rightarrow 0}\sup_{|x-y|<\eta}\!\left(\mathcal{W}_{1}(s,x)\!-\!\mathcal{W}_{2}(t,y)\!-\!\frac{|t-s|^{2}}{2\varepsilon_{2}}\!-\!\frac{\delta}{e^{\rho s}}\!\left(|x|^{m+1}\!+1\right)\!-\!\frac{\delta}{e^{\rho t}}\!\left(|y|^{m+1}\!+1\right)\!-\!\frac{r}{s}\!-\!\frac{r}{t}\right).
\end{align*}
It is easy to see that
\begin{align}\label{eq:Limitm1m2}
\lim_{\varepsilon_{1}\rightarrow 0}m_{1}(\varepsilon_{1},\varepsilon_{2})=m_{2}(\varepsilon_{2}),\quad\lim_{\varepsilon_{2}\rightarrow 0}m_{2}(\varepsilon_{2})=m_{0}.
\end{align}
Note that for $x,y\in\overline{\mathscr{O}}$ with $|x|$ and $|y|$ large enough, $\Psi(t,x,s,y)$ becomes negative. On the other hand, \eqref{eq:MaxPsiGamma0} guarantees that $m_{1}>\gamma_{0}/2$. Hence, $\Psi$ achieves its maximum, which is at least $\gamma_{0}/2$, at some $(\bar{s},\bar{x},\bar{t},\bar{y})$, in certain bounded region. Therefore,
\begin{align*}
m_{1}(\varepsilon_{1},\varepsilon_{2})=\Psi(\bar{t},\bar{x},\bar{s},\bar{y})\leq m_{1}(2\varepsilon_{1},2\varepsilon_{2})-\frac{|\bar{x}-\bar{y}|^{2}}{4\varepsilon_{1}}-\frac{|\bar{t}-\bar{s}|^{2}}{4\varepsilon_{2}},
\end{align*}
and so, by \eqref{eq:Limitm1m2},
\begin{align*}
\lim_{\varepsilon_{2}\rightarrow 0}\,\lim_{\varepsilon_{1}\rightarrow 0}\left(\frac{\left|\bar{x}-\bar{y}\right|^{2}}{4\varepsilon_{1}}+\frac{|\bar{t}-\bar{s}|^{2}}{4\varepsilon_{2}}\right)=0.
\end{align*}

We now show that $\bar{s},\bar{t}\in(0,T)$ and that $\bar{x},\bar{y}\in\mathscr{O}$. From the expression of $\Psi$, it is easy to see that $\bar{s}>0$ and that $\bar{t}>0$. Next, assume that $\bar{s}=T$. By \eqref{eq:PolyAssumpWi}, for fixed $r>0$, $\delta>0$ and $\rho>1$ satisfying \eqref{eq:MaxPsiGamma0}, we can choose $R>0$ such that $|\bar{x}|\leq R$ or $|y|\leq R$, since otherwise $\Psi(\bar{t},\bar{x},\bar{s},\bar{y})$ would achieve a negative value when $|\bar{x}|$ and $|\bar{y}|$ are both large enough, contradicting \eqref{eq:MaxPsiGamma0}. Without loss of generality, we assume that $|\bar{x}|\leq R$. Next, since both $\mathcal{W}_{1}$ and $\mathcal{W}_{2}$ are uniformly continuous on $[0,T]\times(\overline{\mathscr{O}}\cap\{|x|\leq R+1\})$, for any $0<\epsilon<r/(2(R+1)T)$, there exists $\lambda_{0}>0$, such that whenever $(s,x)$ and $(t,y)\in[0,T]\times(\overline{\mathscr{O}}\cap\{|x|\leq R+1\})$ with $|t-s|<\lambda_{0}$ and $|x-y|<\lambda_{0}$,
\begin{align*}
\left|\mathcal{W}_{i}(s,x)-\mathcal{W}_{i}(t,y)\right|<\epsilon,\quad i=1,2.
\end{align*}
Hence, for $\varepsilon_{1}>0$ and $\varepsilon_{2}>0$ small enough so that $|\bar{s}-\bar{t}|=(T-\bar{t})<\lambda_{0}$ and $|\bar{x}-\bar{y}|<\min(\lambda_{0},r/(2(R+1)T),1)$, it follows that $|\bar{y}|\leq R+1$, and setting $\bar{x}=(p_{1},z_{1},\theta_{1})$ and $\bar{y}=(p_{2},z_{2},\theta_{2})$,
\begin{align*}
\left|\mathcal{W}_{1}(\bar{s},\bar{x})-\mathcal{W}_{2}(\bar{t},\bar{y})\right|&\leq\left|\mathcal{W}_{1}(T,\bar{x})-\mathcal{W}_{2}(T,\bar{x})\right|+\left|\mathcal{W}_{2}(\bar{s},\bar{x})-\mathcal{W}_{2}(\bar{t},\bar{y})\right|\\
&\leq|p_{1}z_{1}-p_{2}z_{2}|+\epsilon\leq(\left|\bar{x}\right|+\left|\bar{y}\right|)\left|\bar{x}-\bar{y}\right|+\epsilon\leq\frac{2r}{T}.
\end{align*}
Thus,
\begin{align*}
\mathcal{W}_{1}(\bar{s},\bar{x})-\mathcal{W}_{2}(\bar{t},\bar{y})-\frac{r}{\bar{s}}-\frac{r}{\bar{t}}\leq 0,
\end{align*}
which again contradicts \eqref{eq:MaxPsiGamma0}. Therefore, we must have $\bar{t}<T$. Similarly, $\bar{s}<T$ and also $\bar{x},\bar{y}\in\mathscr{O}$.

\medskip
\noindent
\textbf{Step 2.} We now apply Theorem \ref{thm:IshiiLemma} to obtain some contradiction. Set $\mathcal{Q}_{1}=\mathcal{Q}_{2}=(0,T)\times(\mathscr{O}\cap\{|x|<R+1\})$, $\mathcal{Q}=\mathcal{Q}_{1}\times\mathcal{Q}_{2}$, and define
\begin{align*}
\widetilde{\mathcal{W}}_{1}(s,x)&=\mathcal{W}_{1}(s,x)-\delta e^{-\rho s}\left(|x|^{m+1}+1\right)-\frac{r}{s},\\
\widetilde{\mathcal{W}}_{2}(t,y)&=\mathcal{W}_{2}(t,y)+\delta e^{-\rho t}\left(|y|^{m+1}+1\right)+\frac{r}{t},\\
\varphi(s,x,t,y)&=\frac{|x-y|^{2}}{2\varepsilon_{1}}+\frac{|t-s|^{2}}{2\varepsilon_{2}}.
\end{align*}
The arguments in Step 1 above show that $\widetilde{\mathcal{W}}_{1}-\widetilde{\mathcal{W}}_{2}-\varphi$ achieves a local maximum at
$(\bar{s},\bar{x},\bar{t},\bar{y})\in\mathcal{Q}$. By Theorem \ref{thm:IshiiLemma}, since $\overline{\mathcal{J}}^{2}_{-}\widetilde{\mathcal{W}}_{2}=-\overline{J}^{2}_{+}(-\widetilde{\mathcal{W}}_{2})$, there exist
$\widetilde{G}_{1},\widetilde{G}_{2}\in\mathscr{S}^{4}$ such that
\begin{align*}
\left(\left(
\frac{\bar{s}-\bar{t}}{\varepsilon_{2}},\frac{\bar{x}-\bar{y}}{\varepsilon_{1}}\right)^{T},\widetilde{G}_{1}\right)\in\overline{\mathcal{J}}^{2}_{+}\widetilde{\mathcal{W}}_{1}(\bar{s},\bar{x}),\quad\left(\left(\frac{\bar{s}-\bar{t}}{\varepsilon_{2}},\frac{\bar{x}-\bar{y}}{\varepsilon_{1}}\right)^{T},\widetilde{G}_{2}\right)\in\overline{\mathcal{J}}^{2}_{-}\widetilde{\mathcal{W}}_{2}(\bar{t},\bar{y}),
\end{align*}
and that
\begin{align}\label{eq:IshiiLemmaIneqTildeW1W2}
\left(\begin{array}{cc} \widetilde{G}_{1} & 0 \\ 0 & -\widetilde{G}_{2} \end{array}\right)\leq G_{\varphi}+\varepsilon_{1}G_{\varphi}^{2},
\end{align}
where (setting $z:=(s,x,t,y)$)
\begin{align*}
G_{\varphi}=D_{zz}\varphi(\bar{s},\bar{x},\bar{t},\bar{y})=\left(\begin{array}{cccc} \varepsilon_{2}^{-1} & 0 & -\varepsilon_{2}^{-1} & 0 \\ 0 & \varepsilon_{1}^{-1}I_{3} & 0 & -\varepsilon_{1}^{-1}I_{3} \\ -\varepsilon_{2}^{-1} & 0 & \varepsilon_{2}^{-1} & 0 \\ 0 & -\varepsilon_{1}^{-1}I_{3} & 0 & \varepsilon_{1}^{-1}I_{3} \end{array}\right),
\end{align*}

Taking submatrices by omitting the elements of the first and the fifth rows and columns of the matrices on both sides of \eqref{eq:IshiiLemmaIneqTildeW1W2} leads to
\begin{align}\label{eq:UpperG1G2}
\left(\begin{array}{cc} G_{1} & 0 \\ 0 & -G_{2} \end{array}\right)\leq\frac{1}{\varepsilon_{1}}\left(\begin{array}{cc} I_{3} & -I_{3}\\ -I_{3} & I_{3} \end{array}\right)+2\varepsilon_{1}\left(\begin{array}{cc} \varepsilon_{1}^{-2}I_{3} & -\varepsilon_{1}^{-2}I_{3} \\ -\varepsilon_{1}^{-2}I_{3} & \varepsilon_{1}^{-2}I_{3} \end{array}\right)=\frac{3}{\varepsilon_{1}}\left(\begin{array}{cc} I_{3} & -I_{3} \\ -I_{3} & I_{3} \end{array}\right),
\end{align}
where $G_{1},G_{2}\in\mathscr{S}^{3}$ are submatrices of $\widetilde{G}_{1}$ and $\widetilde{G}_{2}$, respectively, obtained by omitting the first row and the first column. We claim that
\begin{align}\label{eq:barstxyD12pmW}
\left(\frac{\bar{s}-\bar{t}}{\varepsilon_{2}},\frac{\bar{x}-\bar{y}}{\varepsilon_{1}},G_{1}\right)\in\overline{\mathcal{D}}^{(1,2)}_{+}\mathcal{W}_{1}(\bar{s},\bar{x}),\quad\left(\frac{\bar{s}-\bar{t}}{\varepsilon_{2}},\frac{\bar{x}-\bar{y}}{\varepsilon_{1}},G_{2}\right)\in\overline{\mathcal{D}}^{(1,2)}_{-}\mathcal{W}_{2}(\bar{t},\bar{y}).
\end{align}
In fact, by the very definition of $\overline{\mathcal{J}}^{2}_{+}\widetilde{\mathcal{W}}_{1}(\bar{s},\bar{x})$, there exist $(s_{n},s_{n})\in\mathcal{Q}_{1}$ and $(q_{n},r_{n},\widetilde{G}_{n})\in\mathcal{J}^{2}_{+}\widetilde{\mathcal{W}}_{1}(s_{n},x_{n})$,
\begin{align}\label{eq:LimitsxqrGn}
\lim_{n\rightarrow\infty}(s_{n},x_{n})=(\bar{s},\bar{x}),\quad\lim_{n\rightarrow\infty}\left(q_{n},r_{n},\widetilde{G}_{n}\right)=\left(\left(\frac{\bar{x}-\bar{y}}{\varepsilon_{1}},\frac{\bar{s}-\bar{t}}{\varepsilon_{2}}\right)^{T},\widetilde{G}_{1}\right).
\end{align}
Hence for any $n\in\mathbb{N}$, $h\in\mathbb{R}$ and $y\in\mathbb{R}^{3}$,
\begin{align*}
\widetilde{\mathcal{W}}_{1}(s_{n}+h,x_{n}+y)-\widetilde{\mathcal{W}}_{1}(s_{n},x_{n})\leq q_{n}h+r_{n}\cdot y+\frac{1}{2}(h,y^{T})^{T}\cdot\widetilde{G}_{n}(h,y^{T})^{T}+o\left(|h|^{2}+|y|^{2}\right).
\end{align*}
Letting $G_{n}$ be the submatrix of $\widetilde{G}_{n}$ obtained by omitting the first row and the first column, we have
\begin{align*}
\widetilde{\mathcal{W}}_{1}(s_{n}+h,x_{n}+y)-\widetilde{\mathcal{W}}_{1}(s_{n},x_{n})\leq q_{n}h+p_{n}\cdot y+\frac{1}{2}y\cdot G_{n}y+o\left(|h|+|y|^{2}\right),
\end{align*}
and so $(q(n),r(n),G_{n})\in\mathcal{D}^{(1,2)}_{+}\widetilde{\mathcal{W}}_{1}(s_{n},x_{n})$. Together with \eqref{eq:LimitsxqrGn}, this shows the first part of \eqref{eq:barstxyD12pmW}. The second part of \eqref{eq:barstxyD12pmW} can be verified similarly.

Now by the very definitions of $\widetilde{\mathcal{W}}_{1}$ and $\widetilde{\mathcal{W}}_{2}$,
\begin{align*}
\left(\frac{\bar{s}-\bar{t}}{\varepsilon_{2}}+\varphi_{1,s}(\bar{s},\bar{x}),\,\frac{\bar{x}-\bar{y}}{\varepsilon_{1}}+D_{x}\varphi_{1}(\bar{s},\bar{x}),\,G_{1}+D_{xx}\varphi_{1}(\bar{s},\bar{x})\right)&\in\overline{\mathcal{D}}^{(1,2)}_{+}\mathcal{W}_{1}(\bar{s},\bar{x}),\\
\left(\frac{\bar{s}-\bar{t}}{\varepsilon_{2}}-\varphi_{2,t}(\bar{t},\bar{y}),\,\frac{\bar{x}-\bar{y}}{\varepsilon_{1}}-D_{y}\varphi_{2}(\bar{t},\bar{y}),\,G_{2}-D_{yy}\varphi_{2}(\bar{t},\bar{y})\right)&\in\overline{\mathcal{D}}^{(1,2)}_{-}\mathcal{W}_{2}(\bar{t},\bar{y}),
\end{align*}
where $\varphi_{1}(s,x):=\delta e^{-\rho s}(|x|^{m+1}+1)+r/s$ and $\varphi_{2}(t,y):=\delta e^{-\rho t}(|y|^{m+1}+1)+r/t$. It follows from \eqref{eq:SuperDiff} and \eqref{eq:Subdiff} that
\begin{align*}
\frac{\bar{t}-\bar{s}}{\varepsilon_{2}}-\varphi_{1,s}(\bar{s},\bar{x})+\mathcal{H}\left(\bar{s},\bar{x},\frac{\bar{x}-\bar{y}}{\varepsilon_{1}}+D_{x}\varphi_{1}(\bar{s},\bar{x}),G_{1}+D_{xx}\varphi_{1}(\bar{s},\bar{x})\right)&\leq
0,\\
\frac{\bar{t}-\bar{s}}{\varepsilon_{2}}+\varphi_{2,t}(\bar{t},\bar{y})+\mathcal{H}\left(\bar{t},\bar{y},\frac{\bar{x}-\bar{y}}{\varepsilon_{1}}-D_{y}\varphi_{2}(\bar{t},\bar{y}),G_{2}-D_{yy}\varphi_{2}(\bar{t},\bar{y})\right)&\geq
0.
\end{align*}
The above two inequalities immediately lead to
\begin{align}\label{eq:UpperHamilton}
\mathcal{H}\left(\bar{s},\bar{x},\hat{r}_{1},\widehat{G}_{1}\right)-\mathcal{H}\left(\bar{t},\bar{y},\hat{r}_{2},\widehat{G}_{2}\right)+\rho\delta e^{-\rho\bar{s}}\left(|\bar{x}|^{m+1}+1\right)+\rho\delta e^{-\rho\bar{t}}\left(|\bar{y}|^{m+1}+1\right)\leq -\frac{2r}{T^{2}},
\end{align}
where
\begin{align*}
\hat{r}_{1}=\hat{r}_{1}(\bar{s},\bar{x}):=\frac{\bar{x}-\bar{y}}{\varepsilon_{1}}+D_{x}\varphi_{1}(\bar{s},\bar{x}),\quad \hat{r}_{2}=\hat{r}_{2}(\bar{t},\bar{y}):=\frac{\bar{x}-\bar{y}}{\varepsilon_{1}}-D_{y}\varphi_{2}(\bar{t},\bar{y}),\\
\widehat{G}_{1}=\widehat{G}_{1}(\bar{s},\bar{x}):=G_{1}+D_{xx}\varphi_{1}(\bar{s},\bar{x}),\quad\widehat{G}_{2}=\widehat{G}_{2}(\bar{t},\bar{y}):=G_{2}-D_{yy}\varphi_{2}(\bar{t},\bar{y}).
\end{align*}

Now for any $u\in\mathscr{U}$, recalling the notations of $\vec{f}$, $\vec{a}$ and $\mathcal{L}$ in Section \ref{subsec:StoContProb}, we have
\begin{align}
&\left(\!-\vec{f}(\bar{s},\bar{x},u)\!\cdot\!\hat{r}_{1}\!-\!\frac{1}{2}\text{tr}\!\left(\vec{a}(\bar{s},\bar{x},u)\widehat{G}_{1}\!\right)\!-\!\mathcal{L}(\bar{x},u)\!\right)\!-\!\left(\!-\vec{f}(\bar{t},\bar{y},u)\!\cdot\!\hat{r}_{2}\!-\!\frac{1}{2}\text{tr}\!\left(\vec{a}(\bar{t},\bar{y},u)\widehat{G}_{2}\!\right)\!-\!\mathcal{L}(\bar{y},u)\!\right)\nonumber\\
&\quad\geq -\left|\vec{f}(\bar{s},\bar{x},u)\cdot\hat{r}_{1}-\vec{f}(\bar{t},\bar{y},u)\cdot\hat{r}_{2}\right|-\frac{1}{2}\left(\text{tr}\!\left(\vec{a}(\bar{s},\bar{x},u)\widehat{G}_{1}-\vec{a}(\bar{t},\bar{y},u)\widehat{G}_{2}\right)\right)-\left|\mathcal{L}(\bar{x},u)-\mathcal{L}(\bar{y},u)\right|\nonumber\\
\label{eq:DecompHsxHty} &\quad=:-\mathcal{I}_{1}-\mathcal{I}_{2}-\mathcal{I}_{3}.
\end{align}
By Assumption \ref{assump:Coefficients}$-$(i), $\mathcal{I}_{1}$ can be estimated via
\begin{align}
\mathcal{I}_{1}&\leq\left|\left(\vec{f}(\bar{s},\bar{x},u)-\vec{f}(\bar{t},\bar{y},u)\right)\cdot\frac{\bar{x}-\bar{y}}{\varepsilon_{1}}\right|+\left|\vec{f}(\bar{s},\bar{x},u)\cdot D_{x}\varphi_{1}(\bar{s},\bar{x})\right|+\left|\vec{f}(\bar{t},\bar{y},u)\cdot D_{x}\varphi_{2}(\bar{t},\bar{y})\right|\nonumber\\
&\leq\left|\left(\vec{f}(\bar{s},\bar{x},u)-\vec{f}(\bar{s},\bar{y},u)\right)\cdot\frac{\bar{x}-\bar{y}}{\varepsilon_{1}}\right|+\left|\left(\vec{f}(\bar{s},\bar{y},u)-\vec{f}(\bar{t},\bar{y},u)\right)\cdot\frac{\bar{x}-\bar{y}}{\varepsilon_{1}}\right|\nonumber\\
&\quad\,+\delta(m+1)\left|\bar{x}\right|^{m-1}\sup_{u\in\mathscr{U}}\left|\vec{f}(\bar{s},\bar{x},u)\cdot\bar{x}\right|+\delta(m+1)\left|\bar{y}\right|^{m-1}\sup_{u\in\mathscr{U}}\left|\vec{f}(\bar{t},\bar{y},u)\cdot\bar{y}\right|\nonumber\\
\label{eq:CompPrinEstI1} &\leq O\left(\frac{\left|\bar{x}-\bar{y}\right|^{2}}{\varepsilon_{1}}\right)+K_{1}\left(1+|x|^{m+1}+|y|^{m+1}\right),
\end{align}
where $K_{1}>0$ is a constant depending on $\delta$ and $K$ (the Lipschitz constant in Assumption \ref{assump:LipCondbSigma}), but independent of $\rho$. Next, the last term $\mathcal{I}_{3}$ in \eqref{eq:DecompHsxHty} can be estimated via
\begin{align}\label{eq:CompPrinEstI3}
\mathcal{I}_{3}\leq\kappa\mathscr{L}^{\gamma}\left|\bar{x}-\bar{y}\right|.
\end{align}
Finally, for the second term $\mathcal{I}_{2}$, first by \eqref{eq:UpperG1G2},
\begin{align*}
\text{tr}\left(\vec{a}(\bar{s},\bar{x},u)G_{1}-\vec{a}(\bar{t},\bar{y},u)G_{2}\right)\leq\frac{3}{\varepsilon_{1}}\text{tr}\left(\left(\vec{\sigma}(\bar{s},\bar{x},u)-\vec{\sigma}(\bar{t},\bar{y},u)\right)\left(\vec{\sigma}^{T}(\bar{s},\bar{x},u)-\vec{\sigma}^{T}(\bar{t},\bar{y},u)\right)\right),
\end{align*}
and together with Assumption \ref{assump:Coefficients}$-$(i) as well as the uniform continuity of $\vec{\sigma}$ in $[0,T]\times(\overline{\mathscr{O}}\cap\{|x|\leq R+1\})$, this leads to
\begin{align}
\mathcal{I}_{2}&\leq\frac{3}{2\varepsilon_{1}}\left(\text{tr}\left(\left(\vec{\sigma}(\bar{s},\bar{x},u)-\vec{\sigma}(\bar{t},\bar{y},u)\right)\left(\vec{\sigma}^{T}(\bar{s},\bar{x},u)-\vec{\sigma}^{T}(\bar{t},\bar{y},u)\right)\right)\right)\nonumber\\
&\quad\,+\frac{1}{2}\left(\text{tr}\left(\vec{a}(\bar{s},\bar{x},u)D_{xx}\varphi_{1}(\bar{s},\bar{x})+\vec{a}(\bar{t},\bar{y},u)D_{xx}\varphi_{2}(\bar{t},\bar{y})\right)\right)\nonumber\\
\label{eq:CompPrinEstI2} &\leq
O\left(\frac{\left|\bar{x}-\bar{y}\right|^{2}}{\varepsilon_{1}}\right)+K_{2}\left(1+|\bar{x}|^{m+1}+|\bar{y}|^{m+1}\right),
\end{align}
where $K_{2}>0$ is a constant depending on $\delta$ and $K$ (the Lipschitz constant in Assumption \ref{assump:Coefficients}$-$(i)), but independent of $\rho$. Combining \eqref{eq:DecompHsxHty}$-$\eqref{eq:CompPrinEstI2}, we obtain (denoting $K_{0}:=K_{1}+K_{2}$)
\begin{align*}
\mathcal{H}\!\left(\bar{s},\bar{x},\hat{r}_{1},\widehat{G}_{1}\right)-\mathcal{H}\!\left(\bar{t},\bar{y},\hat{r}_{2},\widehat{G}_{2}\right)\geq -K_{0}\left(1+\left|\bar{x}\right|^{m+1}+\left|\bar{y}\right|^{m+1}\right)+O\left(\left|\bar{x}-\bar{y}\right|\right)+O\left(\frac{\left|\bar{x}-\bar{y}\right|^{2}}{\varepsilon_{1}}\right),
\end{align*}
which, together with \eqref{eq:UpperHamilton}, leads to
\begin{align*}
-K_{0}\!\left(1\!+\!\left|\bar{x}\right|^{m+1}\!+\!\left|\bar{y}\right|^{m+1}\right)+O\!\left(\left|\bar{x}\!-\!\bar{y}\right|\right)+O\!\left(\frac{\left|\bar{x}\!-\!\bar{y}\right|^{2}}{\varepsilon_{1}}\right)+\frac{\rho\delta}{e^{\rho\bar{s}}}\!\left(\left|\bar{x}\right|^{m+1}\!\!+\!1\right)+\frac{\rho\delta}{e^{\rho\bar{t}}}\!\left(\left|\bar{y}\right|^{m+1}\!\!+\!1\right)\leq -\frac{2r}{T^{2}}.
\end{align*}
Choose $\rho>1$ large enough so that
\begin{align*}
\rho\delta e^{-\rho\bar{s}}\left(\left|\bar{x}\right|^{m+1}+1\right)+\rho\delta e^{-\rho\bar{t}}\left(\left|\bar{y}\right|^{m+1}+1\right)-K_{0}\left(1+\left|\bar{x}\right|^{m+1}+\left|\bar{y}\right|^{m+1}\right)>0.
\end{align*}
By taking $\varepsilon_{2}\rightarrow 0$ and then $\varepsilon_{1}\rightarrow 0$, we finally obtain that
\begin{align*}
0=\limsup_{\varepsilon_{1}\rightarrow 0}\left[\limsup_{\varepsilon_{2}\rightarrow 0}\left(O\left(\left|\bar{x}-\bar{y}\right|\right)+O\left(\frac{\left|\bar{x}-\bar{y}\right|^{2}}{\varepsilon_{1}}\right)\right)\right]\leq-\frac{2r}{T^{2}}<0,
\end{align*}
which is clearly a contradiction. The proof is now complete.\hfill $\Box$

\end{document}